\UseRawInputEncoding 

\documentclass{article}

\usepackage{color,graphicx,times}
\usepackage{amssymb,amscd}
\usepackage{pb-diagram}     
\usepackage{amsmath}
\usepackage{amsthm}
\usepackage{graphicx}

\makeatother

\let \ttorg \tt \def \tt{\ttorg \obeyspaces}

\begin{document}

\newcommand{\Across}{\raisebox{-0.25\height}{\includegraphics[width=0.5cm]{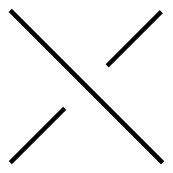}}}
\newcommand{\Bcross}{\raisebox{-0.25\height}{\includegraphics[width=0.5cm]{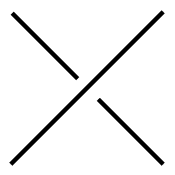}}}
\newcommand{\Asmooth}{\raisebox{-0.25\height}{\includegraphics[width=0.5cm]{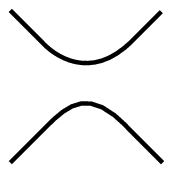}}}
\newcommand{\Bsmooth}{\raisebox{-0.25\height}{\includegraphics[width=0.5cm]{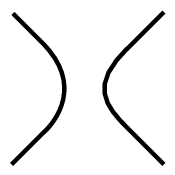}}}
\newcommand{\Rcurl}{\raisebox{-0.25\height}{\includegraphics[width=0.5cm]{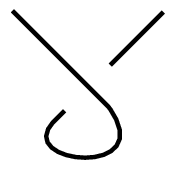}}}
\newcommand{\Lcurl}{\raisebox{-0.25\height}{\includegraphics[width=0.5cm]{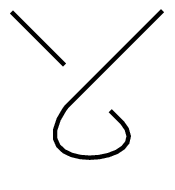}}}
\newcommand{\Arc}{\raisebox{-0.25\height}{\includegraphics[width=0.5cm]{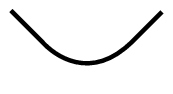}}}

\title{\Large \bf Combinatorial Knot Theory and the Jones Polynomial}
\author{Louis H. Kauffman\\ Department of Mathematics, Statistics \\ and Computer Science (m/c
249)    \\ 851 South Morgan Street   \\ University of Illinois at Chicago\\
Chicago, Illinois 60607-7045\\ $<$kauffman@uic.edu$>$\\}

\maketitle

\thispagestyle{empty}

\noindent {\bf Abstract.} This paper is  an introduction to combinatorial knot theory via state summation models for the Jones polynomial and its generalizations. 
It is also a story about the developments that ensued in relation to the discovery of the Jones polynomial and a remembrance of Vaughan Jones and his mathematics. \\



\section{Introduction}
This paper describes developments in knot theory
that were inspired by the Conway Skein Theory and the Jones polynomial. These developments involve a wide range of fields and ideas and provide an opportunity to see mathematics, physics and natural science through a remarkable window.

The paper consists in five sections beyond the introduction. In Section 2  we recall the skein theory of John Horton Conway and how this led to the author's discovery of a state summation model for the Alexander-Conway polynomial.
This state summation was later instrumental in finding a state summation model for the Jones polynomial. This section discusses the 
Homflypt and Kauffman two variable polynomials and the role of the connection structure version of the Temperley-Lieb algebra that the author discovered along with the bracket model for the Jones polynomial.
These early state summation models were, as we now know, the tip of an iceberg. Section 2 discusses Penrose tensor networks and their role in construction of so-called quantum link invariants via solutions
to the Yang-Baxter equation.  Section 3 describes Witten's breakthrough, giving a model of the Jones polynomial via functional integration and quantum field theory. This section outlines how the Witten approach is related to Vassiliev invariants (defined in this section) and how
the Vassiliev invariants can be used to trace a direct line from combinatorial knot theory to Lie algebra (supporting the weight systems of Vassiliev invariants). The quantum field theoretic approach to link invariants
stands in the middle between deformed Lie algebras (quantum groups and Hopf algebras) and purely combinatorial approaches using just Lie algebra alone. Needless to say, we cannot tell all the detail in this section, but we do show how one arrives at the diagrammatic footprint of the Jacobi identity by finding the relations (the four-term relations) implied for Vassiliev evaluations by the invariance under the Reidemeister moves.
Here the knot theory gives a hint about categories of diagrams that underlie both algebra and topology. Section 4 introduces virtual knot theory, a theory of knots embedded in thickened surfaces, and discusses one
invariant, the Arrow polynomial, that is a non-trivial generalization of the Jones polynomial giving subtle and interesting information about virtual knots and links. The Arrow polynomial is of particular interest in this story, since it is a generalization of the Jones polynomial that has extra strength for knot and links embedded in thickened surfaces. This is a first step toward Jones polynomials for knots and links in three-manifolds. Section 5 discusses Khovanov homology. The Khovanov Category of a knot or link has already been introduced in Section 2 in terms of the states of the Kauffman bracket. We then address the question: How can one extract topological information about knots and links from 
this category of states. Toward an answer we discuss the Cobordism Category construction of Dror Bar-Natan and show how it leads to factoring surface cobordisms by the Four Tube Relation. This gives a diagrammatic and categorical structure that is the backbone of Khovanov homology. We describe how the Bar-Natan Cobordism Category yields a natural proof of the invariance of the Khovanov homology under Reidemeister moves. It is a proof that lines up with the original proof of the invariance of the bracket polynomial. We end with the remark that the Heegaard-Floer Homology of Oszvath and Szabo is (by their work) combinatorially rooted in the Formal Knot Theory states with which we began this essay in describing our exploration of the Alexander polynomial.  

Section 6 is a printing of a letter from Vaughan Jones to the author, written in the fall of 1986. The letter is full of insight and valuable to this day.

\begin{figure}[t]
	\centering
	\includegraphics[width=.5\textwidth]{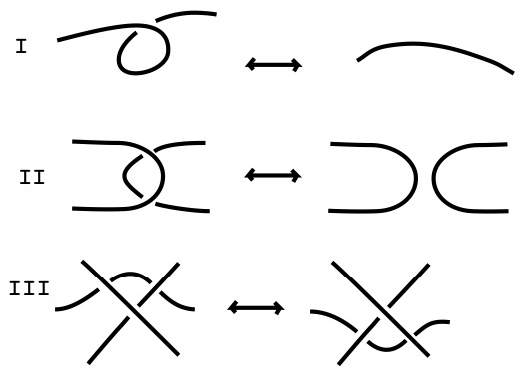}
	\caption{Reidemeister moves.}
	\label{moves}
\end{figure}

\begin{figure}[t]
	\centering
	\includegraphics[width=.5\textwidth]{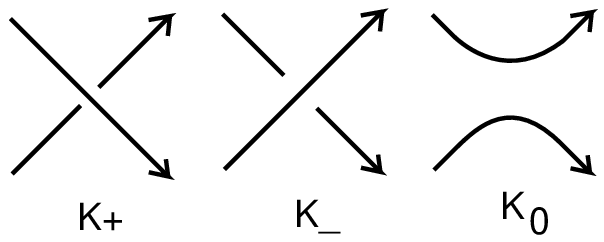}
	\caption{Skein triple.}
	\label{skein}
\end{figure}

\section{Bracket Polynomial and Jones Polynomial}
Before beginning to describe the bracket polynomial and the Jones polynomial, we remark that these invariants of knots and links are based in a diagrammatic approach that 
was discovered by J. W. Alexander and Garland Baird Briggs (1926)  \cite{AlexBriggs}  and by Kurt Reidemeister (1927)  \cite{RM}. These researchers discovered a set of combinatorial moves on diagrams for knots and links such that two links embedded in three dimensional space are ambient isotopic (equivalent by a continuous family of embeddings) if and only if any two projection diagrams of these links are equivalent by the Reidemeister moves. (It has been customary to refer to the moves as Reidemeister moves because they are the foundation of the book ``Knotentheorie" \cite{Reidemeister} by Reidemeister, published in 1934.) See Figure~\ref{moves} for an illustration of each move type. The Reidemeister moves provide a complete planar combinatorial translation of the problems of knot and link equivalence in three dimensional Euclidean space.\\

Vaughan Jones' discovery \cite{JO} of a new Laurent polynomial invariant of oriented knots and links $V_{K}(t)$ came as a bolt out of the blue in 1984. The new invariant was
derived from a representation of the Artin braid group to the Temperley-Lieb algebra, an algebra that had originally been found as a matrix algebra by Temperley and Lieb in
the early 1970's \cite{Baxter}. Jones rediscovered this algebra in studying a construction for von Neumann algebras that produced a tower of algebras from an inclusion $M \subset N$ of von Neumann algebras. Each element in the tower has a projection to its predecessor, and there comes forth an algebra generated by projections $e_1,e_2,e_3,\cdots$ so that
$e_{i}^{2} = e_{i}, e_{i}e_{i \pm1}e_{i} = \tau e_{i}$ ($\tau$ a scalar) and $e_{i}e{_j} = e_{j}e{_i}$ when $|i-j|>1.$ These are the relations for the Temperley-Lieb algebra.  Jones noted the similarity of these relations to the generating relations for the Artin braid group, and he proceeded to find a representation of the Artin braid group to this algebra. There was more.
His work on von Neumann algebras led him to construct a trace function on this Temperley-Lieb algebra (let it be understood that a trace function $tr$ satisfies $tr(ab) = tr(ba)$ for products of elements $a$ and $b$ in the algebra) and to wonder how this would interact with the braid group representation. Consultation with Joan Birman \cite{Margalit} led Jones to construct a knot invariant from that trace by using the Markov Theorem (telling when a trace on the braid group can yield a knot invariant) and the Jones polynomial $V_{K}(t)$ was born. Jones discovered that $V_{K}(t)$ can often tell the topological difference between knots and their mirror images. The Alexander polynomial cannot distinguish mirror images,  and so the new polynomial was not the 
Alexander polynomial \cite{Alex}.  Furthermore, the new invariant was related to von Neumann algebras and to statistical mechanics. Jones was very generous with his speculations and results about the new polynomial and its context, and gave many talks on it during its first year in the world.

John Horton Conway \cite{Conway} had reformulated the structure of the  Alexander polynomial $\Delta_{K}(t)$ \cite{Alex} published by James W. Alexander in 1928. The Conway version $\nabla_{K}(z)$ is determined by the skein axioms:
\begin{enumerate}
	\item $\nabla_{K}(z) = \nabla_{K'}(z)$ whenever $K$ and $K'$ are ambient isotopic oriented links.
	\item $\nabla_{K}(z) = 1$ if $K$ is an unknotted single loop.
	\item $\nabla_{K_{+}} - \nabla_{K_{-}} = z \nabla_{K_{0}}$ whenever  $K_{+}, K_{-}, K_{0}$ are three diagrams that differ only at one local site where in one there is a positive crossing, in the next there is a negative crossing and in the third there are parallel arcs, as shown in Figure~\ref{skein}.
\end{enumerate}
The Conway polynomial can be computed from just the diagram of the knot or link $K$ by simple recursion and no other mathematical apparatus. Its relation with the classical Alexander polynomial is encapsulated by the formula $\nabla_{K}(\sqrt{t} - 1/\sqrt{t}) \doteq \Delta_{K}(t)$ where $\doteq$ denotes equality up to a factor of $\pm t^{N}$ for some integer $N$.

Jones showed that his new polynomial satisfied a skein relation similar to the Conway skein relation. He proved that 
$$t^{-1} V_{K_{+}}(t) - t V_{K_{-}}(t) = (\sqrt{t} - 1/\sqrt{t}) V_{K_{0}}(t).$$ This led a group of people (Hoste, Ocneanu, Millett, Freyd, Lickorish, Yetter, Przytycki and Trawczk) \cite{homfly,pt} to, 
independently and in pairs, define an invariant two-variable common generalization $P_{K}(a,z)$ of the Jones and Alexander--Conway polynomials, with the elegant skein relation 
$$aP_{K_{+}}(a,z) - a^{-1}P_{K_{-}}(a,z) = z P_{K_{0}}(a,z).$$ It is known by its acronym as the Homflypt Polynomial. Subsequently I discovered another two variable polynomial that uses
unoriented diagrams, sometimes called the Kauffman Polynomial \cite{KauffmanPoly}, but we are getting ahead of the story.

\begin{figure}[t]
	\centering
	\includegraphics[width=\textwidth]{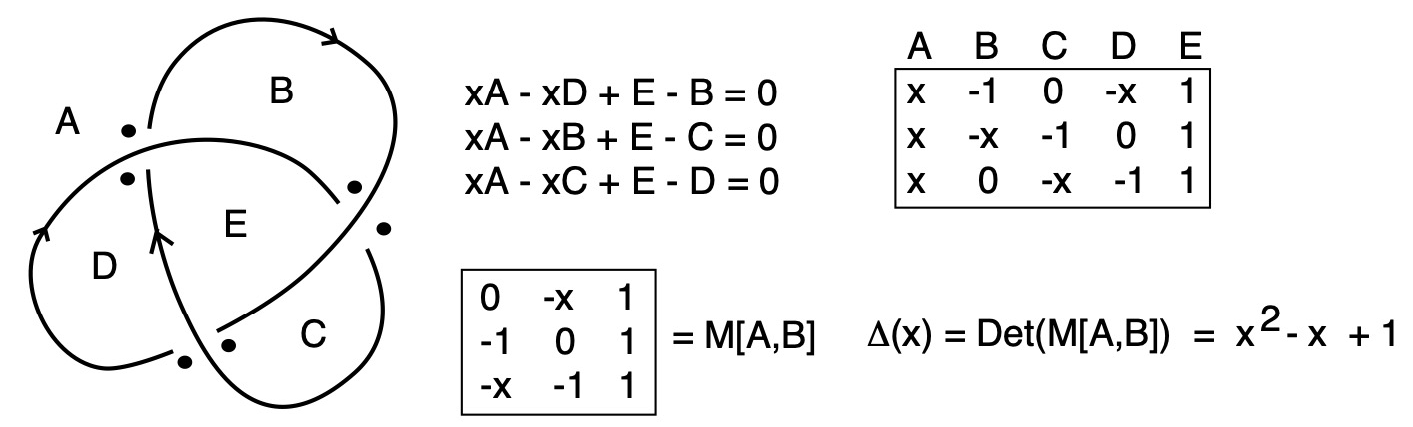}
	\caption{Alexander's algorithm for the Alexander polynomial.}
	\label{alex}
\end{figure}

\begin{figure}[t]
	\centering
	\includegraphics[width=4.5cm]{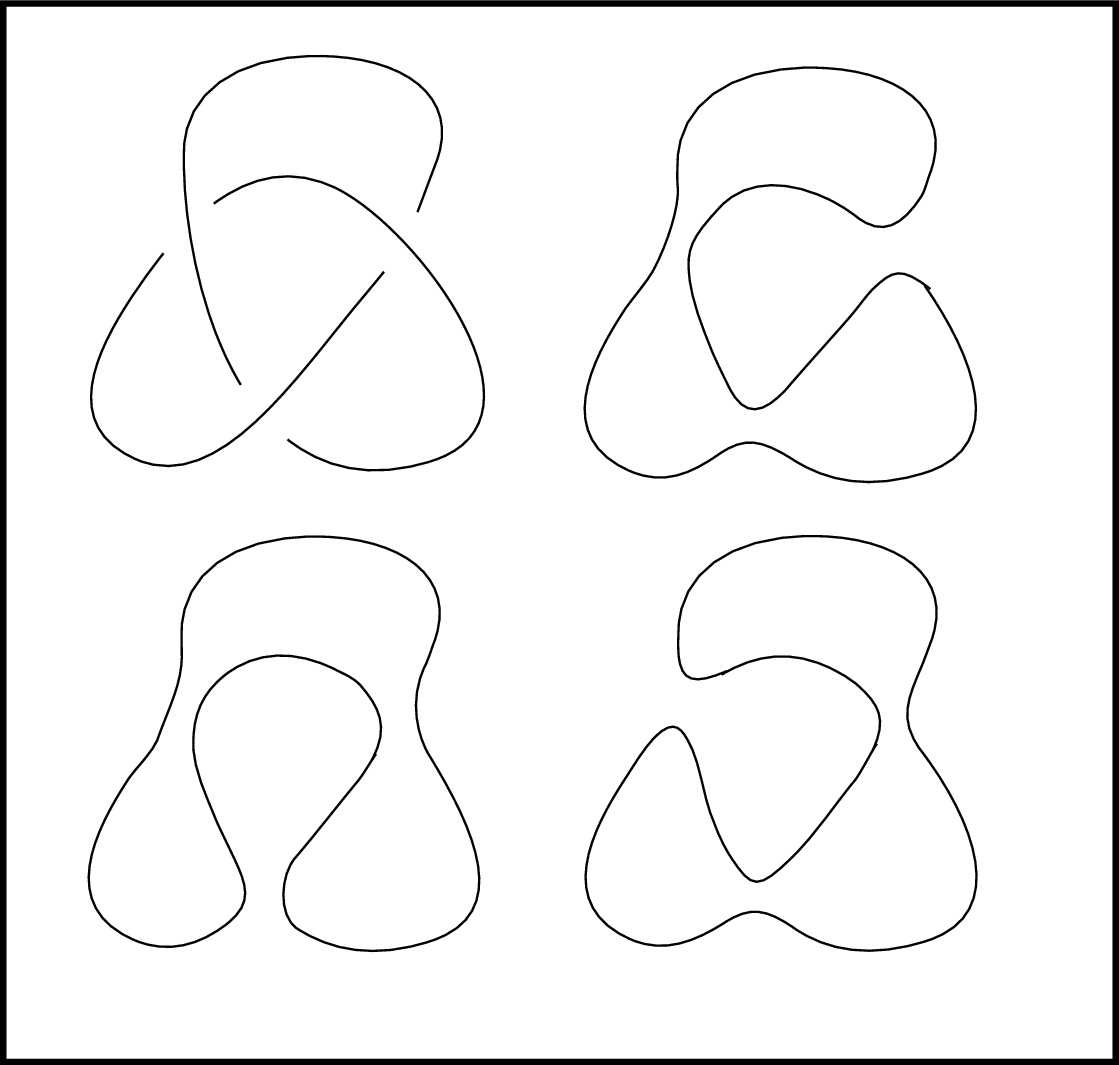}
	\caption{Jordan-Euler trails for the trefoil diagram.}
	\label{loops}
\end{figure}

Since 1980 I had been working on what I called a {\it state summation model} for the Alexander-Conway polynomial. This was a combinatorial formula for the Alexander-Conway polynomial, that was a sum over all the ways to smooth the crossings in the knot diagram (unoriented) so that one obtains a single Jordan curve in the plane as a result.
I call these {\it Jordan-Euler Trails} on the knot diagram (Euler because it was Euler who first considered walks on a graph that use each edge once).
See Figure~\ref{loops}. Each Jordan-Euler trail contributes a term to the polynomial. This summation over combinatorics produces the invariant polynomial and is analogous to the sum over the states of a physical system called its {\it partition function} \cite{Baxter, KA89}. 

\begin{figure}[t]
	\centering
	\includegraphics[width=.5\textwidth]{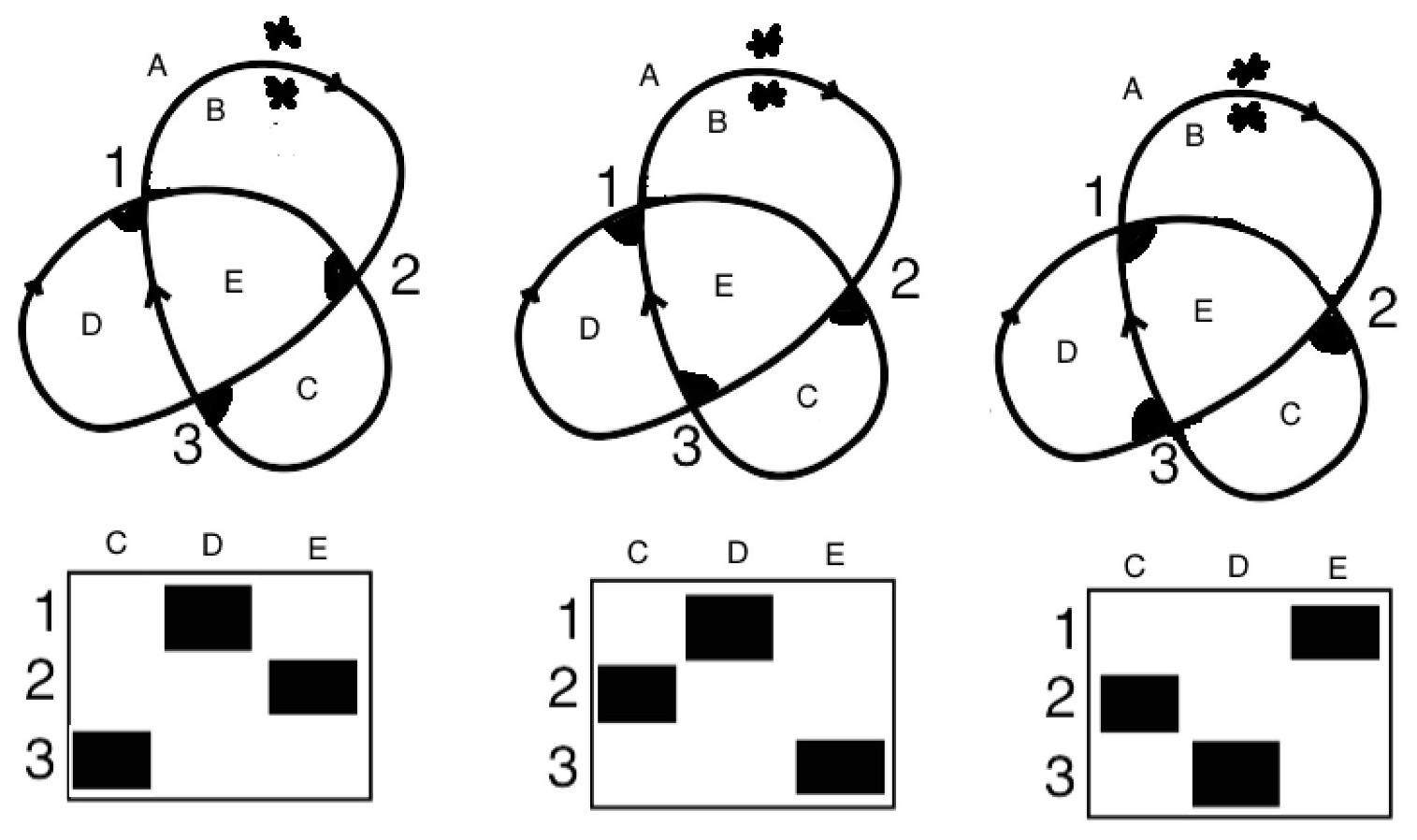}
	\caption{States for the Alexander algorithm.}
	\label{states}
\end{figure}

\begin{figure}[t]
	\centering
	\includegraphics[width=.2\textwidth]{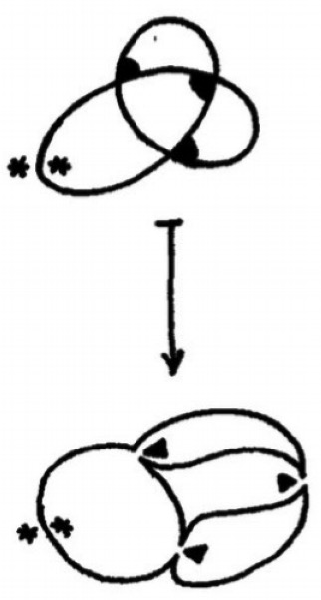}
	\caption{States and trails.}
	\label{trail}
\end{figure}

The idea that the Alexander-Conway polynomial should come from this combinatorics came from understanding the remarkable structure in Alexander's original paper. In Figure~\ref{alex} I illustrate Alexander's original algorithm for his polynomial. As the reader can glean from the figure, a module is generated by the regions of the diagram and there is a relation among the regions incident to a crossing. These relations assemble into a matrix and the determinant of this matrix with two columns deleted (that correspond to two adjacent regions in the diagram) gives the polynomial. The Alexander polynomial is determined up to a sign and a power of its variable $t,$ and is invariant under the Reidemeister moves. Figure~\ref{states} illustrates my reformulation of Alexander's determinant as a state summation. The states consist in choices, by the regions, of crossings in the diagram. These choices are indicated by the black triangular markers on the knot diagram. Each state corresponds to a term in the expansion of the determinant. A term in the expansion of the determinant corresponds having each column choose a unique row of the matrix. The figure shows how the markers in each state correspond to positions in the (Alexander) matrix whose product gives a term in the determinant.  To see this look at the square black markers in the figure that indicate positions in the three by three matrix.
For example, just below the first state on the left there is a three by three matrix with a black mark at the 1D position. This corresponds to the state marker at node 1 for region D. The reader will need to examine the figure in the light of these remarks to see more. The terms in the determinant expansion correspond to these marker states, and the marker states correspond to the Jordan-Euler trails in Figure~\ref{loops} that we have already mentioned. The relationship between marker states and trails is shown in Figure~\ref{trail} where each marker is used to smooth a crossing and the trail appears from this smoothing process. The figure shows what is meant by a marker corresponding to a smoothing. Other features come into play, not the least of which is that the permutation signs in the determinant expansion can be obtained directly from a parity in the state diagrams. In this way, a state summation emerges that has a life of its own and which can be used to support the structure of the Alexander-Conway polynomial and is a conceptual tool for investigating its properties.

I succeeded in finding the state summation model and had published a book about it \cite{FKT}. When I heard about the Homflypt polynomial, I was sure that there must be a state summation of that invariant, generalizing the one I had found for Alexander-Conway. And I searched, to no avail, for such a model in the fall of 1984 and the spring of 1985.
Then, late in the summer of 1985, Lickorish, Millett and Ho \cite{HO} discovered a one variable skein polynomial based on unoriented diagrams, and I  realized that it could be generalized to a two-variable skein polynomial invariant by using a framing variable making the polynomial an invariant of regular isotopy (the equivalence relation generated by the second two Reidemeister moves). This new polynomial invariant satisfies the skein relation
$$L_{\Across} + L_{\Bcross} = z(L_{\Asmooth} + L_{\Bsmooth})$$ 
coupled with the behaviour under curls
$$L_{\Rcurl} = a L_{\Arc},$$
$$L_{\Lcurl} = a^{-1} L_{\Arc}.$$

A few days after this discovery, I was on a plane to Italy to visit  Massimo Ferri in Bologna. On the plane, it occurred to me that I could look for a state summation
specialization of this polynomial in the form that,  given an unoriented link diagram $K$, there is 
associated to it a well-defined Laurent polynomial in the variable $A$,  $\langle K \rangle (A).$
$$\langle \Across \rangle=A \langle \Asmooth \rangle + B\langle \Bsmooth \rangle$$
$$\langle K \, \bigcirc \rangle= d \langle K \rangle, $$
$$\langle \bigcirc \rangle = 1.$$
Here we take $B = A^{-1}$ and  $d = -A^{2} -A^{-2}.$
The last equation guarantees that the bracket evaluates to unity on an unknotted circle. The choice of  specialization of $B$ and $d$ guarantees that the bracket is invariant
under the second and third Reidemeister moves, as we explain below.
The small diagrams indicate parts of otherwise identical larger knot or link diagrams. The two types of smoothing (local diagram with no crossing) in this formula are said to be of type $A$ ($A$ as above) and type $B$ ($B$ as above).

\begin{figure}[t]
	\centering
	\includegraphics[width=6cm]{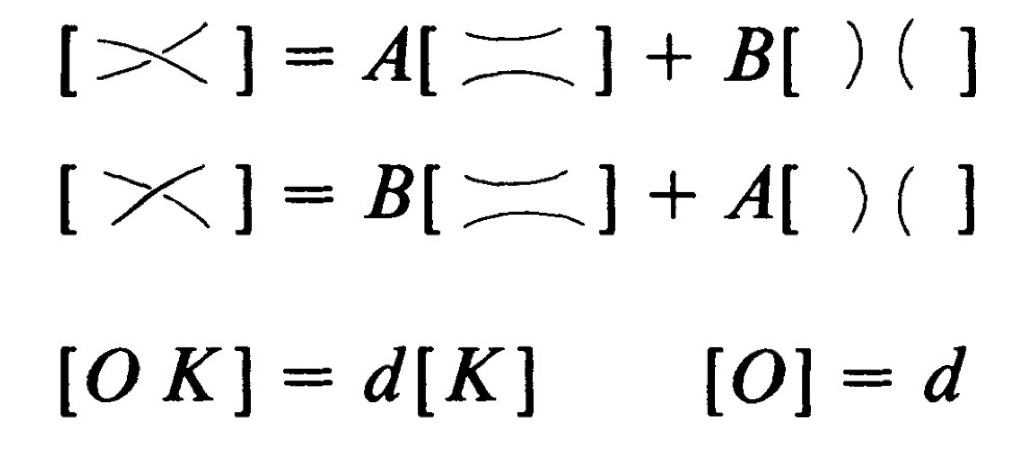}
	\caption{Bracket axioms.}
	\label{axioms}
\end{figure}

\begin{figure}[t]
	\centering
	\includegraphics[width=\textwidth]{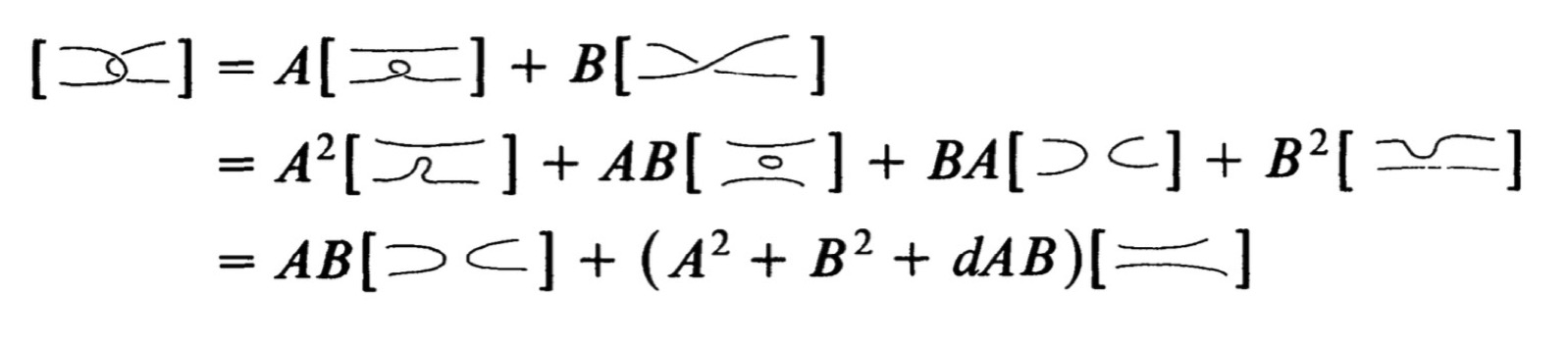}
	\caption{Bracket expansion for second Reidemeister move.}
	\label{R2}
\end{figure}

\begin{figure}[t]
	\centering
	\includegraphics[width=.7\textwidth]{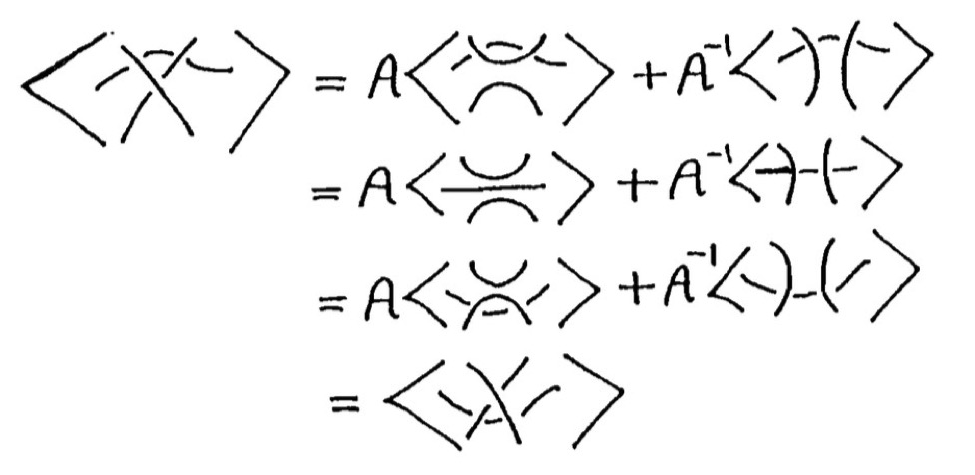}
	\caption{Bracket expansion for third Reidemeister move, given invariance under second Reidemeister move.}
	\label{R3}
\end{figure}

Figure~\ref{axioms} shows the axioms for the bracket and Figure~\ref{R2} shows the expansion of the bracket for the form of the Reidemeister two move. As the reader can see, the bracket will be invariant
under the move when $B = A^{-1}$ and  $d = -A^{2} -A^{-2}.$ We see from Figure~\ref{R3} that once $A,B$ and $d$ are chosen so that we have invariance under the second Reidemeister move, it 
follows that there is also invariance under the third Reidemeister move. This is accomplished by expanding on one crossing in the triangular pattern of the move, and then applying invariance under
the second move. Without the specialization of the variables the bracket polynomial can be used as a combinatorial polynomial associated with knot and link diagrams and it is then directly related to
the dichromatic polynomial in graph theory and the partition function of the Potts model in statistical mechanics \cite{KaB,KAM,KA89,KP}.

It is a consequence of this setup that the bracket behaves as below under the curls that are eliminated by the first Reidemeister move.
$$\langle \Rcurl \rangle=(-A^{3}) \langle \Arc \rangle, $$
$$\langle \Lcurl \rangle=(-A^{-3}) \langle \Arc \rangle .$$
Note that $\langle \Across \rangle + \langle \Bcross \rangle = (A+B)(\langle \Asmooth \rangle + \langle \Bsmooth \rangle),$ and in this way the bracket becomes a special case of the 
$L$-polynomial.

The bracket is  often called the {\it Kauffman bracket polynomial} and the $L$-polynomial is called the {\it (two variable) Kauffman polynomial.}
On the plane to Italy, I did not yet realize that the bracket polynomial was an un-normalized version of the Jones polynomial. The plane landed and Massimo sent me off to Venice to 
get touristic experience at once so we could settle down to mathematics as soon as possible. And so it was in Venice that I had the pleasure of realizing that the bracket 
yielded a simple model for the Jones polynomial.
In the normalized version we define $$f_{K}(A) = (-A^{3})^{-wr(K)} \langle K \rangle $$
where the writhe $wr(K)$ is the sum of the oriented crossing signs for a choice of orientation of the link $K.$   One then has that $f_{K}(A)$ is invariant under the Reidemeister moves
(see \cite{KaB,KAM,KA89}) and the original Jones polynomial $V_{K}(t)$ is given by the formula $$V_{K}(t) = f_{K}(t^{-1/4}).$$ The Jones polynomial
has been of great interest since its discovery in 1984 due to its relationships with statistical mechanics, due to its ability to often detect the
difference between a knot and its mirror image and due to the many open problems and relationships of this invariant with other aspects of low
dimensional topology. It was a remarkable experience to realize that it could be defined so simply in terms of the bracket state summation, and that this meant that {\it the Jones polynomial itself takes the form 
	of a partition function in statistical mechanics}, and that it is a knot theoretic relative of the Tutte polynomial (The Tutte polynomial is equivalent to the  the dichromatic polynomial \cite{KA89,KAM}.) with the two smoothings of the knot diagram corresponding to the deletion and contraction of graphical edges.

\medskip
\noindent {\bf The State Summation.} In order to obtain a closed formula for the bracket, we now describe it as a state summation.
Let $K$ be any unoriented link diagram. Define a {\em state} $S$ of $K$ to be the collection of planar loops resulting from  a choice of
smoothing for each  crossing of $K.$ There are two choices ($A$ and $B$) for smoothing a given  crossing, and
thus there are $2^{c(K)}$ states of a diagram with $c(K)$ crossings.
In a state we label each smoothing with $A$ or $A^{-1}$ according to the convention
indicated by the expansion formula for the bracket. These labels are the  {\em vertex weights} of the state.
There are two evaluations related to a state. The first is the product of the vertex weights,
denoted $\langle K|S \rangle .$
The second is the number of loops in the state $S$, denoted  $||S||$.

\begin{figure}[t]
	\centering
	\includegraphics[width=.5\textwidth]{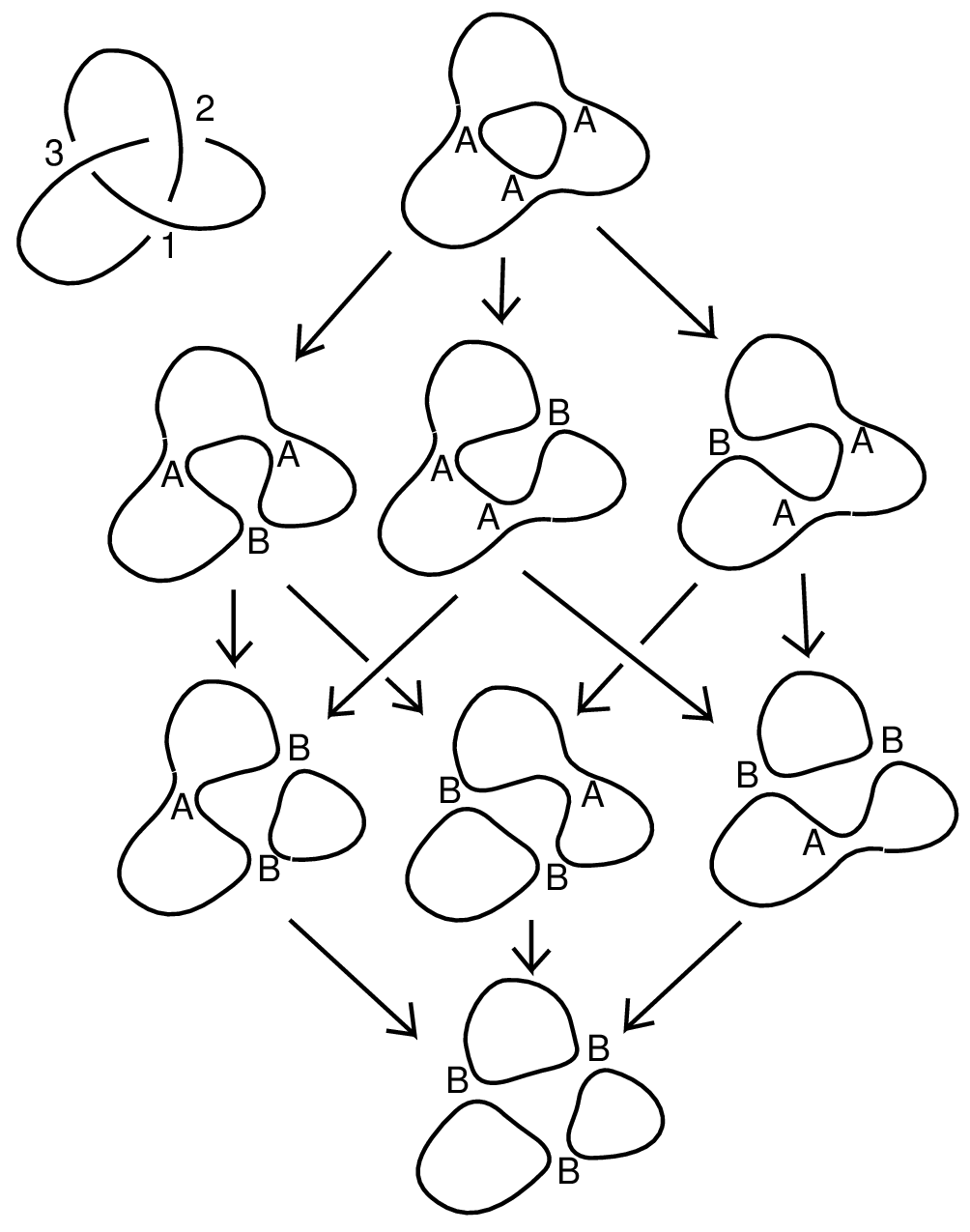}
	\caption{Bracket states and Khovanov complex.}
	\label{cube}
\end{figure}

\noindent 
Define the {\em state summation}, $\langle K \rangle $, by the formula
$$\langle K \rangle  \, = \sum_{S} <K|S> d^{||S||-1}$$
where $d = -A^{2} - A^{-2}.$
This is the state expansion of the bracket. 
In Figure~\ref{cube} we show all the states for the right-handed trefoil knot, labelling the sites
with $A$ or $B$ where $B$ denotes a smoothing that would receive $A^{-1}$ in the state expansion.
Note that in the state enumeration in Figure~\ref{cube} we have organized the states in tiers so that the state
that has only $A$-smoothings is at the top and the state that has only $B$-smoothings is at the bottom. 

\begin{figure}[t]
	\centering
	\includegraphics[width=.5\textwidth]{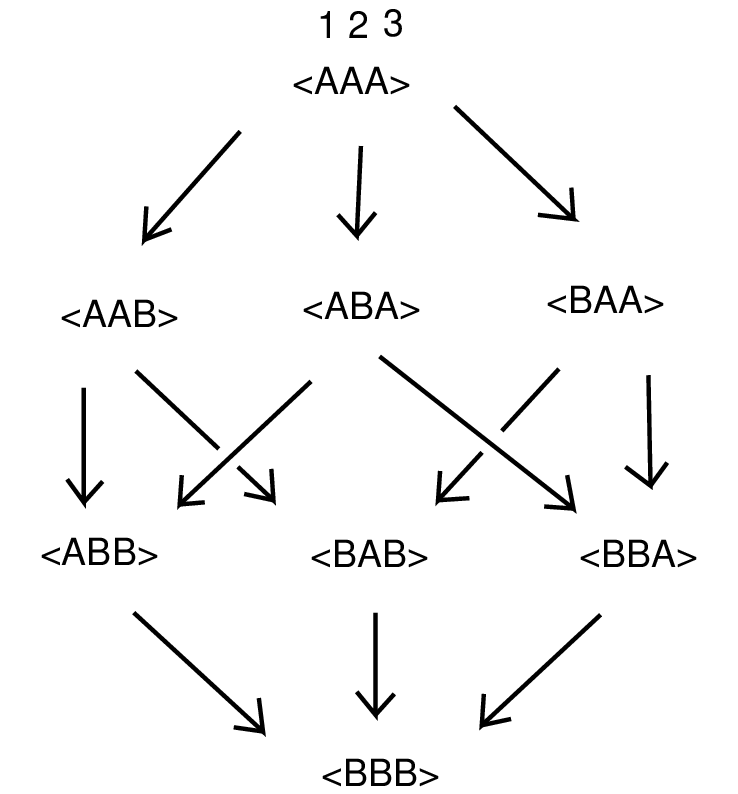}
	\caption{The cube category.}
	\label{cubecat}
\end{figure}

This organization, with arrows taking a state $S$ to a state $S'$ so that $S'$ has one more $B$-smoothing, gives the states the structure of a category. The arrows between the states generate, by composition of arrows, the arrows in the category. The objects are the states and there is an unwritten identity arrow from each object to itself.
Note that an arrow in this figure shows a state changing to another state by resmoothing one $A$-smoothing to a $B$-smoothing. This is the Khovanov Category \cite{Kho} and is the beginning of a breakthrough into Link Homology that occurred in Mikhail Khovanov's work in 1999. In Figure~\ref{cubecat} we illustrate the {\it cube category} that is the framework of this Khovanov category for a knot or link diagram. In the cube category each node of the cube graph is an object in the category and each directed edge is a generating morphism. Here you see a 3-cube, but if the diagram has $n$ crossings, then it will have an $n$-cube category in its background. 

The cube category itself is an example of {\it categorification,} a term for opening up a mathematical structure by turning it into a category. This means that one is respecting certain distinctions that were formerly ignored.
In this case what was ignored is the possibility to order the states by having arrows from $A's$ to $B's!$ So simple, but a new world arises in the production and analysis of the resulting category. You can get a feel for
this sort of movement by thinking of how the algebra of $(A+B)^3 = A^3 + B^3 + 3A^2 B + 3A B^2$ is related to the structure of a cube with side-length $A+B$ where there will be  smaller cubes and parallelepipeds of volume $A^3$, $B^3$, $A^2 B$ and $A B^2.$ The usual algebra does not include the way that these pieces are glued together to form the larger cube.

The cube category of Figure~\ref{cubecat} comes from making categorical sense of the algebraic expression $(A \longrightarrow B)^{3}$ and you can see by looking at the figure how the category does indeed describe the decomposition of the cube into its component cubes and rectangles. Consider that one could make an algebra of expressions like $(A \longrightarrow B)^{2}$ and write
$$(A \longrightarrow B)^{2} = (A \longrightarrow B)(A \longrightarrow B)$$ 
$$= (A \longrightarrow B)A \longrightarrow (A \longrightarrow B)B$$
$$= (AA \longrightarrow BA) \longrightarrow (AB \longrightarrow BB),$$
letting products and arrows distribute across the arrows. In the last diagram we see that the expression $(A \longrightarrow B)^{2}$ has expanded into a {\it higher category}. That is it has an arrow that points between two arrows.
In the higher category the arrows $(AA \longrightarrow BA)$ and $(AB \longrightarrow BB)$ are also objects in the category and there can be an arrow between them. If we raise $(A \longrightarrow B)$ to the third 
power and distribute there will be arrows of higher order still. But we can {\it flatten} an arrow between arrows to form an ordinary category by shifting the higher arrow to ordinary arrows between the end objects. If we do that for our example $$(AA \longrightarrow BA) \longrightarrow (AB \longrightarrow BB),$$ we obtain the $2$-cube category as shown below.
$$
\begin{CD}
AA&@>>>BA\\
@VVV&@VVV\\
AB&@>>>BB
\end{CD}
$$\\

\begin{figure}[t]
	\centering
	\includegraphics[width=.9\textwidth]{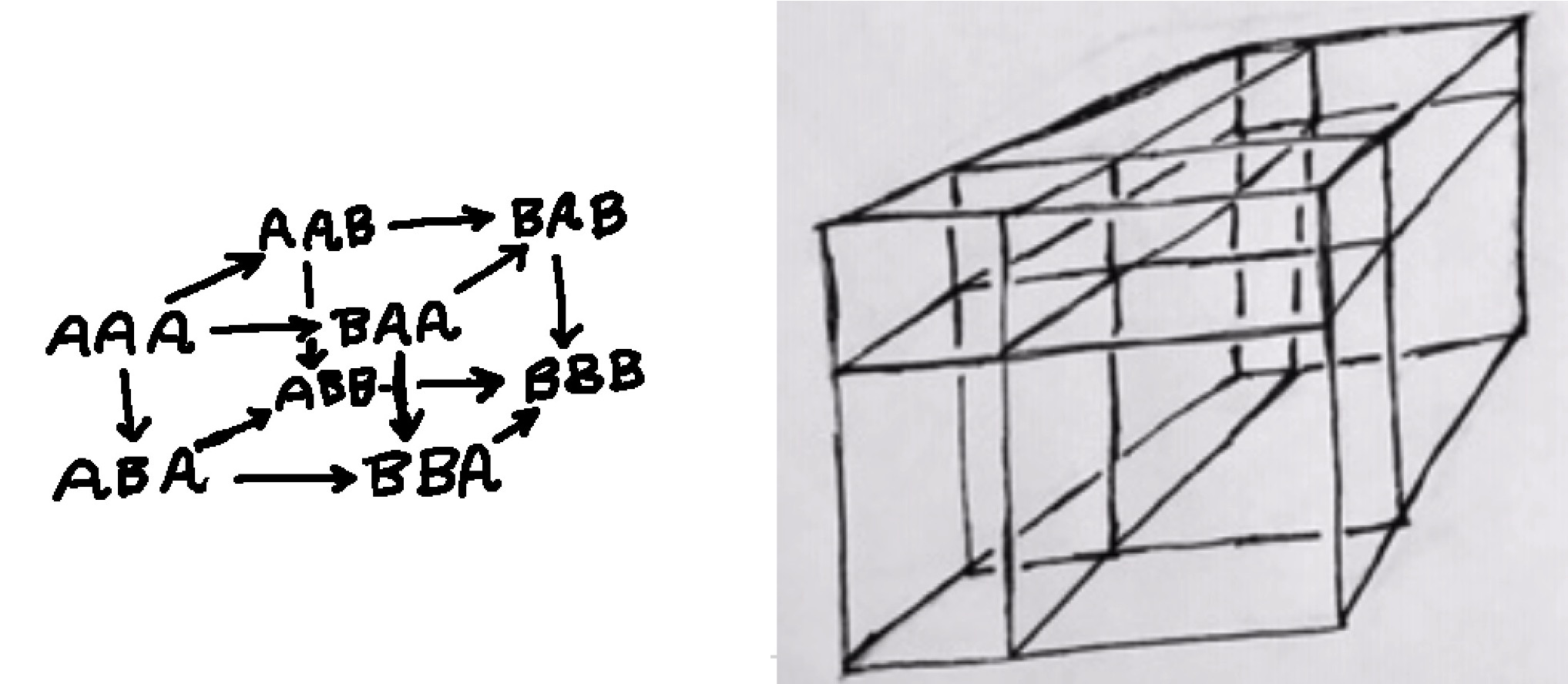}
	\caption{The cube category and the cube.}
	\label{cubecatcube}
\end{figure}

Thus we obtain the $n$-cube category by the prescription ${\cal F} [(A \longrightarrow B)^{n}]$ where ${\cal F}$ denotes the operation of flattening a higher category to a standard category. The language of 
higher categories and categorical algebra lets us describe the detailed decomposition of an $n$-dimensional cube.
One can think of this notion of categorification as a retrograde motion. In the $1500's$, before our modern algebra, mathematicians like Cardano and Tartaglia  had to refer directly to 
the decomposition of a three dimensional cube to enable their solution to the roots of a cubic equation. Categorification embraces old and new mathematical structures in wider patterns. Figure~\ref{cubecatcube} illustrates the cube category juxtaposed with the architecture of a three dimensional cube.
This picture of the Khovanov category gets ahead of our story, and it shows how the state formulation of the Jones polynomial became a seed for future developments. We will return to the Khovanov homology in Section 5.

The states of the bracket comprise all possible smoothings of the diagram and so include the states that were used \cite{FKT} to make a model for the Alexander polynomial.
In fact one can model the Jones polynomial with this restricted set of states. The result is more technical but deeply related to the Tutte polynomial in graph theory \cite{KAM}.
The bracket expansion identity is a knot diagrammatic version of the contraction-deletion relation so central to much of graph theoretic analysis. In this sense the bracket polynomial was a breakthrough between graph theory and knot theory, a breakthrough that is continuing to expand in the present time.

\begin{figure}[t]
	\centering
	\includegraphics[width=.5\textwidth]{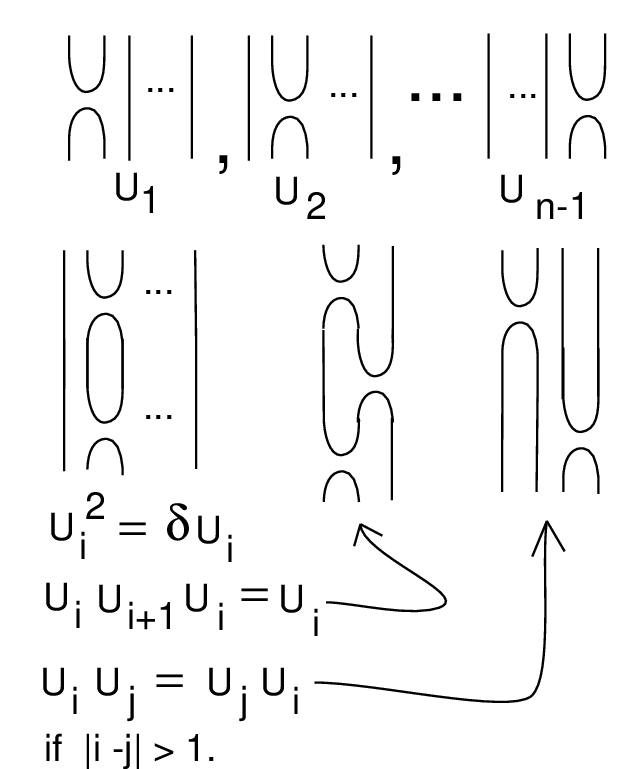}
	\caption{Diagrams for Temperley-Lieb algebra.}
	\label{TL}
\end{figure}

\begin{figure}[t]
	\centering
	\includegraphics[width=.7\textwidth]{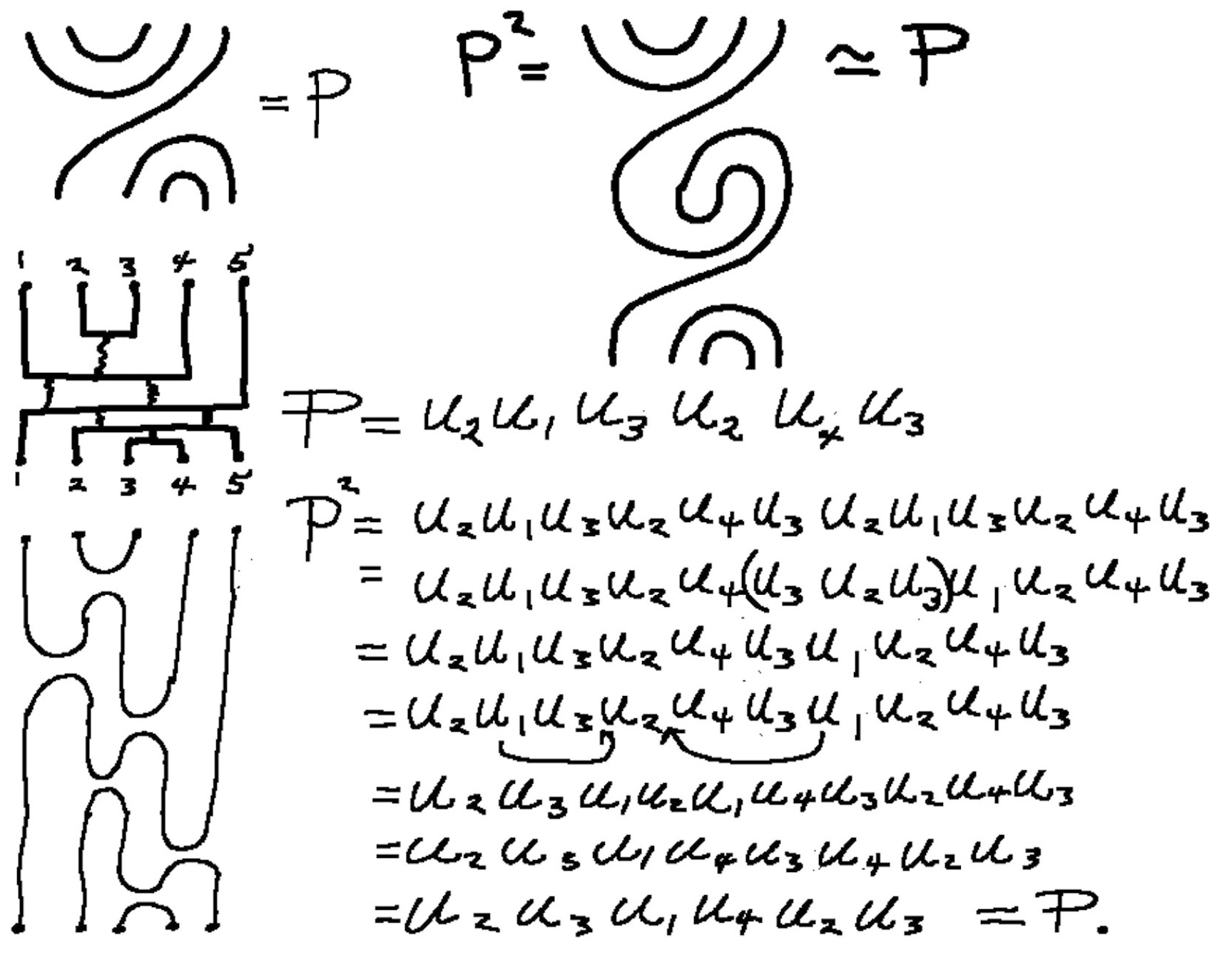}
	\caption{Connection algebra and Temperley-Lieb algebra.}
	\label{connect}
\end{figure}

My story of the fall of 1985 continues. I talked about these discoveries in Bologna and then continued on to Torino where I visited the physicist Mario Rasetti. There, continuing to 
lecture on this material, I discovered that the bracket polynomial was {\it directly} related to the Temperley-Lieb algebra and the other ideas that were involved in the Jones definition of the polynomial. If you apply the bracket formula to a braid you are led to consider some very suggestive diagrams as shown in Figure~\ref{TL}. We obtain a diagrammatic/combinatorial definition of the Temperley-LIeb algebra with the relations in the form $U_{i}^{2} = d U_{i}, U_{i}U_{i\pm1}U_{i} = U_{i}, $  $U_{i}U_{j}=U_{j}U_{i}$ when $|i-j|>1.$ In this form the Temperley-Lieb algebra is a {\it planar connection algebra} with multiplicative generators corresponding to
connections between two rows of points (points on a given row can be connected to one another) under the constraint that the set of connecting arcs embeds in a planar rectangle between the rows. See Figure~\ref{connect}. In this figure we show an example of such a connection structure and how it can be canonically associated with a product of the algebra generators $U_{i}.$ The method for producing the canonical product is to draw the connections in minimal rectangular form and then decorate this form with pairings that will become maximal and minimal in the columns in between the points at the top and the bottom of the diagram. This description will become clear if the reader will view the figure. In that figure we illustrate a connection structure $P$ and show directly that $P^2 = P$ by topological deformation and we show how to translate $P$ into a product of the generators of the Temperley-Lieb algebra and then show that $P^2 = P$ by using the algebraic relations.  A more intensive examination of this relationship shows that the connection algebra is described faithfully by these generators and relations.  See  \cite{KD,TL}.

\begin{figure}[t]
	\centering
	\includegraphics[width=.6\textwidth]{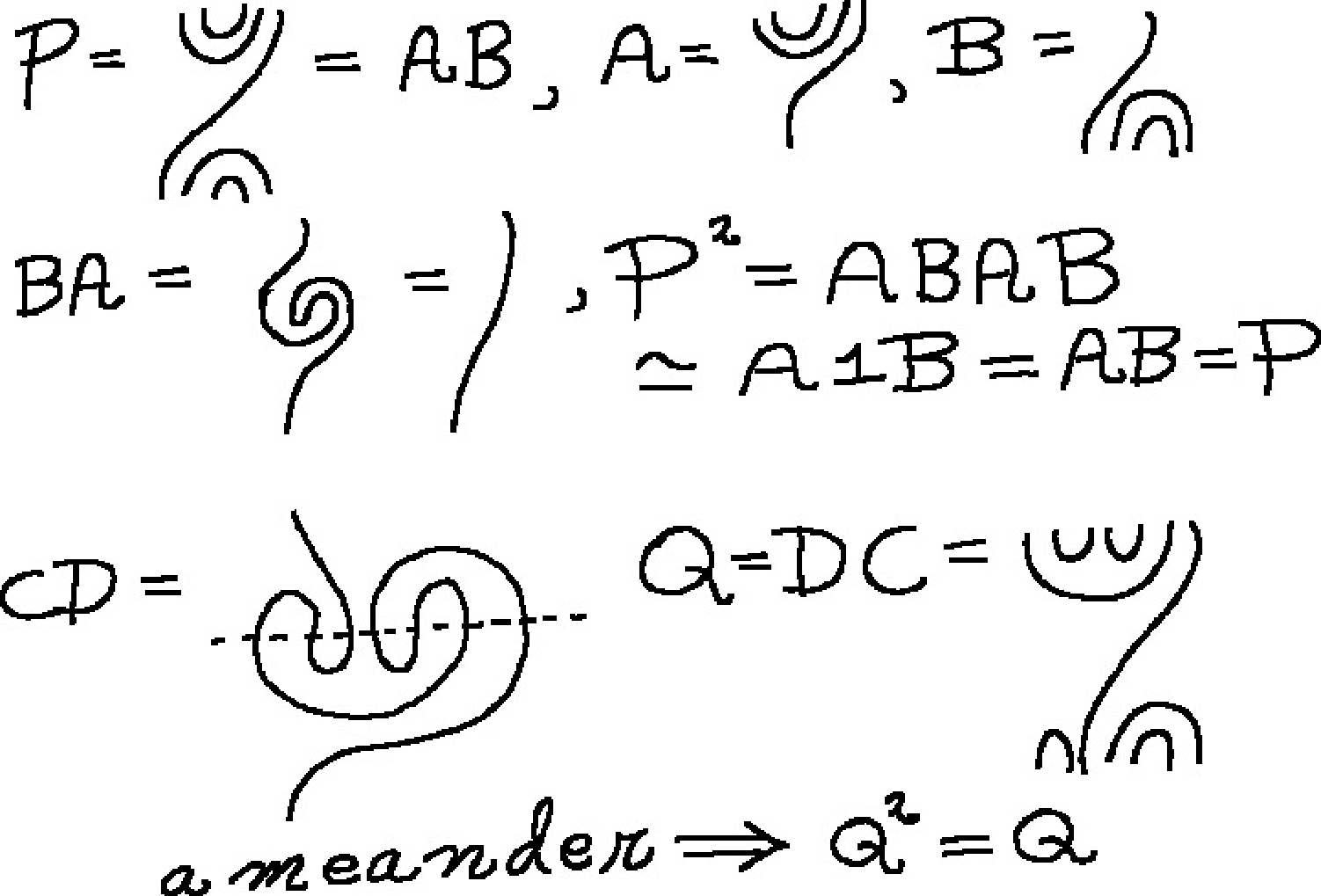}
	\caption{Meanders, projectors and the Temperley-Lieb category.}
	\label{catconnect}
\end{figure}

In Figure~\ref{catconnect} we illustrate how the {\it Temperley-Lieb Connection Category} can illuminate the structure of this last example. In this figure we factor $P = BA$ where $P$ is as in the previous figure and
$B$ and $A$ are morphisms in the Connection Category. Such morphisms are planar connections between two rows of points, but there are different numbers of points in the two rows. In the case of the morphisms
$A$ and $B$ in this figure, the top row of $A$ has five points and the bottom row has one point, while the bottom row of $B$ has five points and the top row of $B$ has one point. Thus we can form the connections
$AB$ and $BA.$ As is apparent from the figure, $AB$ is a morphism from one point to one point and it is the same (topologically) as the identity morphism. But $BA = P$ our previous morphism from five points to five points. We see that the factorization $AB$ is a {\it meander} in the sense that it is the result of drawing a connected curve in the plane and then cutting it with a horizontal line. Classifying meanders is a venerable and 
fascinating combinatorial subject \cite{M,M1}. Here, we see how to make elements of the Temperley-Lieb algebra that are idempotent by using meanders. We start with a meander $M$ and slice it to obtain a factorization of the identity
$1 = CD$. We let $Q= DC$ and we find that $Q^2 = QQ = DCDC = D1C = DC = Q.$ This algebra occurs in the Connection Category and it applies to the Temperley-Lieb algebra by rewriting $Q$ as a product of the 
Temperley-Lieb generators. To find and classify all idempotents of this type in the Temperly-Lieb algebra we have expanded our view to the Connection Category and availed ourselves of the concept of meanders in that realm. This example shows how taking a wider and categorical view can shed light on a question that might be intractable in its original formulation.

\begin{figure}[t]
	\centering
	\includegraphics[width=.5\textwidth]{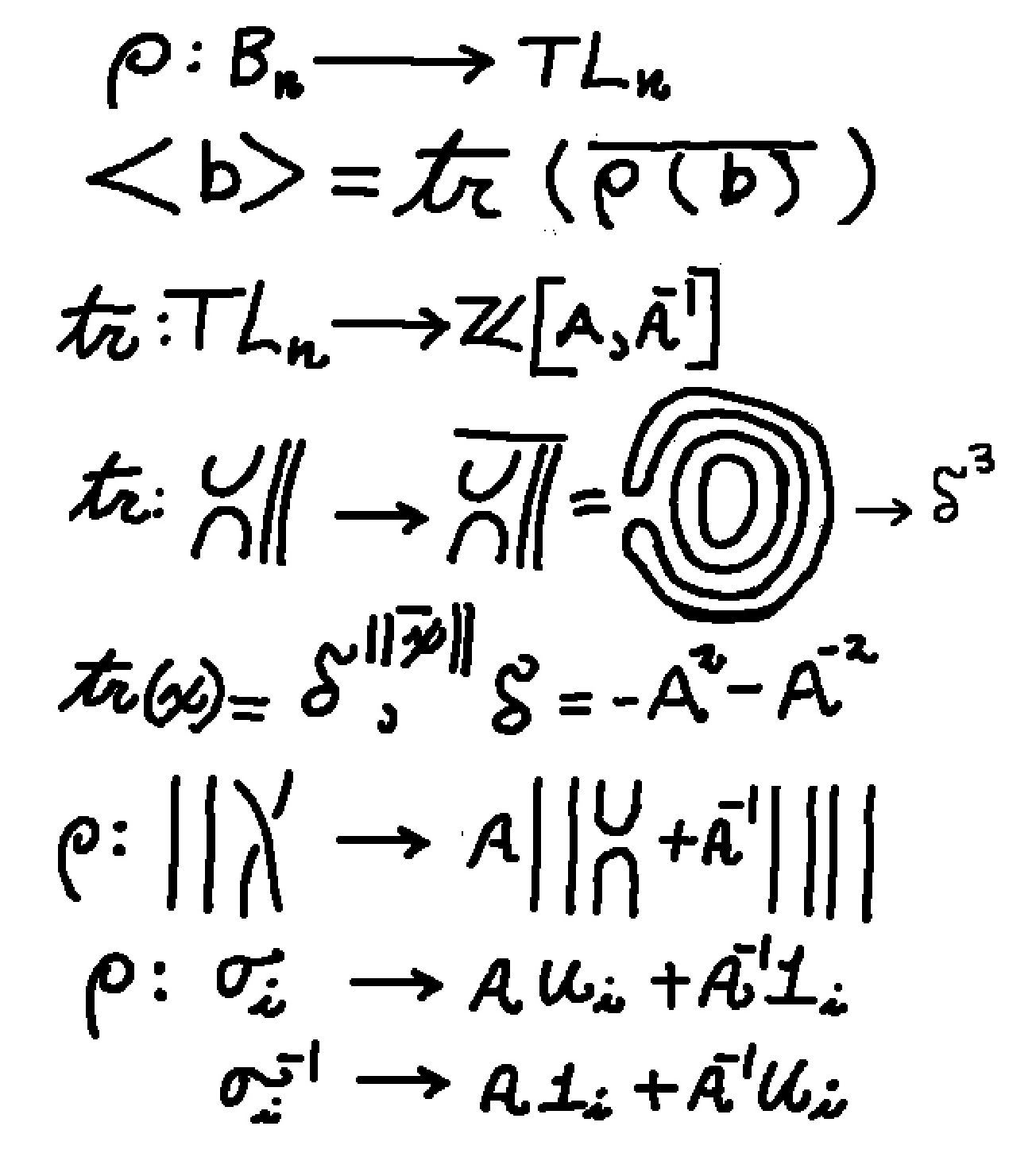}
	\caption{Bracket polynomial via trace on connection algebra.}
	\label{bracketconnect}
\end{figure}

The bracket expansion for a braid  can be regarded as a representation $\rho: B_{n} \longrightarrow TL_{n}$ from the Artin braid group to the Temperley-Lieb algebra, expressed in this mode of planar diagrammatic algebra. Our reformulation of the Jones trace on the abstract Temperley-Lieb algebra corresponds to raising $d$ to he number of loops in the closure of the diagrammatic terms of $\rho(b)$ where $b$ is in $B_{n}.$  Thus the discovery of the bracket model expanded the context of original Jones polynomial. Furthermore, the relationship with statistical mechanics and the Potts model is direct with the bracket model \cite{KaB,KAM,KA89,KP}. The Potts model for planar graphs can be expressed in terms of the bracket formalism and the original relation with the Temperley-Lieb algebra reappears exactly through this combinatorics.

What happened next was an explosion of new mathematics. Jones, Turaev and Reshetikhin and Akutsu and Wadati  discovered many more state summation models and new knot invariants by using 
solutions to the Yang-Baxter equation and formulating all of this in terms of quantum groups and Hopf algebras \cite{AW1,AW2,JO1,JO2,LOMI,RT,Turaev,TuraevViro,Kirby-Melvin}. Then Witten \cite{Witten,Atiyah} discovered a quantum field theoretic interpretation of 
the Jones polynomial and its relatives. Witten's work partially solves the question of a three dimensional topological interpretation of the Jones polynomial. The qualification is that the functional integrals in the Witten approach exist in a physical level of rigor. Much more comes after this.  But I want to end this part of the essay with a quote from a letter that I received 
from Vaughan in October  1986, two years prior to Witten's revolution. It is poignant to see the depth of his intuition
 for this connection of physics, combinatorics, algebra and topology. The entire letter is printed in Section 6.\\
 
\subsection{Letter of Vaughan Jones, October 1986.}

\noindent ``Institut des Hautes Etudes Scientifiques\\
3 Oct, 1986\\
Dear Lou,\\
Since I'm about to talk about it today, I thought I should let you know of a states model for the 2-variable polynomial. ...
The model is very suggestive of the ``real" meaning of the polynomials. $L$ [the diagram] should be replaced by a link in 3-space,
the `states' by functions from $L$ to an $(n+1)$ dimensional Hilbert space ... and the sum over contributing states by an integral with respect to some Wiener measure of an 
interaction term depending on the link in $R^{3}.$ Thus is is an object of {\it gauge quantum field theory on $L$}, the gauge group in this case being $SU(n+1).$
I am morally sure that if one expresses the gauge group by $SO(n+1)$ one will obtain the Kauffman polynomial. And there should be other polynomials for all the 
Coxeter Dynkin diagrams...   One last word - the relation with the fundamental group seems rather suggestive but puzzling at this stage. Converting the `vertex model' described above to an `IRF' mdoel on the planar dual, we see that the states assign numbers to the generators in the Dehn presentation of the fundamental group.This suggests a
relationship I have long suspected between $V$ and representations of $\pi_{1}(S^3 - L)$ into $SU(2)$ tying up hopefully with Casson's invariant.\\
Best wishes,\\
Vaughan"

\subsection{Tensor Networks, Partition Functions and Knot Invariants}
\begin{figure}[t]
	\centering
	\includegraphics[width=.4\textwidth]{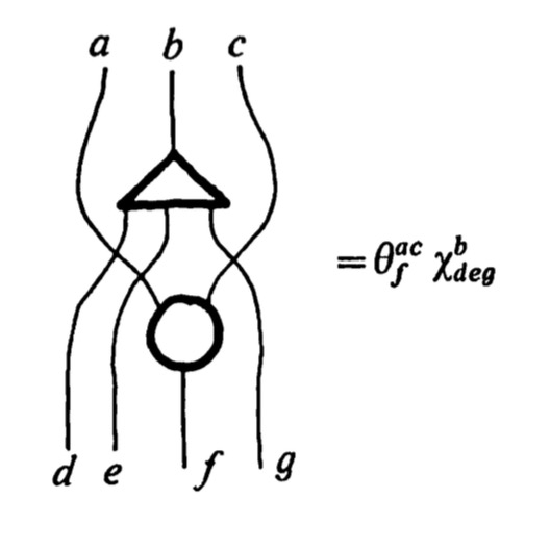}
	\caption{Tensor net.}
	\label{Tens}
\end{figure}

\begin{figure}[t]
	\centering
	\includegraphics[width=.7\textwidth]{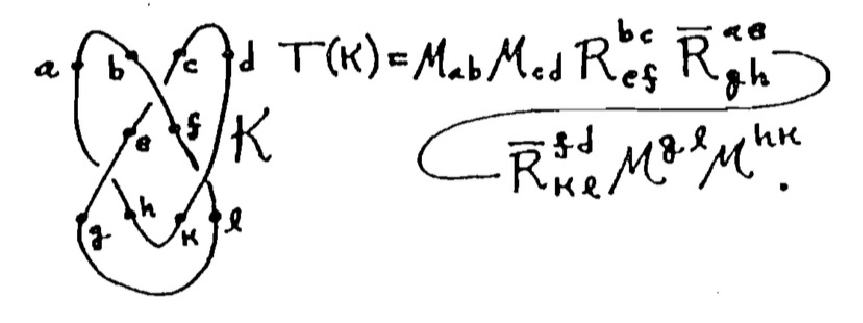}
	\caption{Topological tensor nets.}
	\label{KnotTens}
\end{figure}

A tensor network is a graph $G$ with tensors or matrices associated with each of its nodes and an index set $I$ that can be used to label the edges of the graph. A contraction of the tensor net $G$ is obtained by
assigning fixed indices to all external edges of $G$ and then summing over all possible assignments of indices to internal edges the products of the corresponding matrix entries for the nodes of the graph. See Figure~\ref{Tens} for illustrations of abstract tensor networks. It is often useful to choose a form for the nodes of the graph that is mnemonic for particular uses. For example, in Figure~\ref{KnotTens} we illustrate a tensor network that
is associated with a knot diagram. The basic ingredients in this network are the cups, caps and crossings shown in the figure. Appropriate matrix choices including nodes that satisfy the Yang-Baxter equation, can be used to form the so-called quantum invariants of knots and links.

Tensor networks of the type used here were originally devised by Roger Penrose \cite{Penrose} to handle spin-networks for (quantum) angular momentum, pre-geometry and graph coloring problems. At this stage in the process, I became interested in expressing the state sums that used solutions to the Yang-Baxter equations in terms of tensor networks
\cite{KAT}. {\it Tensor network} is the present term for these structures. In the Penrose work these networks are called {\it abstract tensors}.

\begin{figure}[t]
	\centering
	\includegraphics[width=.3\textwidth]{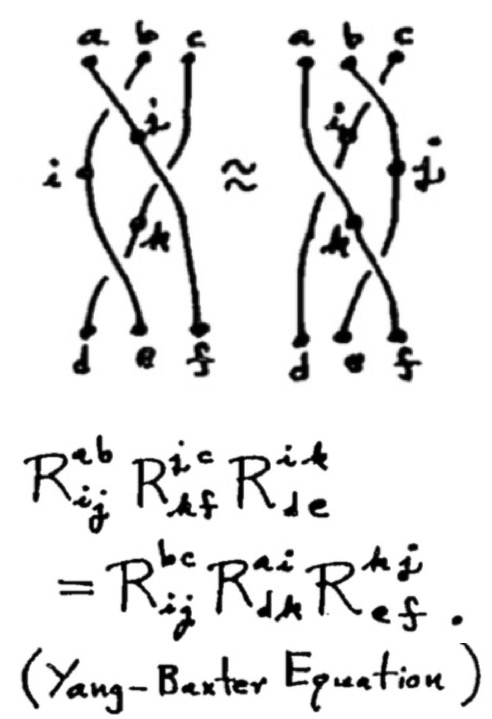}
	\caption{The Yang-Baxter equation.}
	\label{ybe}
\end{figure}

State sum models can be formulated as tensor network contractions on appropriately labeled knot and link diagrams. This point of view includes both the statistical mechanical
viewpoint of Baxter and the quantum field theoretic viewpoint of Yang. In both cases the matrix $R^{ab}_{cd}$ representing a crossing as in Figure~\ref{KnotTens} must satisfy the Yang-Baxter Equation, as shown in Figure~\ref{ybe}. Remarkably,
both Yang and Baxter discovered this equation for reasons that would simplify calculations in statistical mechanics and in quantum field theory. It was only around 1986 that mathematicians began to fully
appreciate that this equation was an expression for the basic braiding relation or so-called third Reidemeister move in knot theory. In this way, statistical mechanics at the hands of Baxter and quantum field theory
at the hands of Yang prepared the basis for an explosion of invariants of knots and links at the hands of topologists and interested physicists.

\begin{figure}[t]
	\centering
	\includegraphics[width=.5\textwidth]{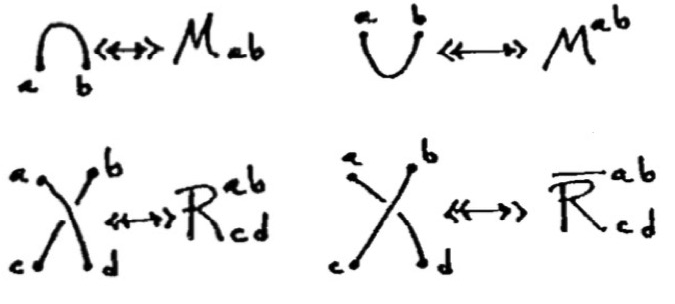}
	\caption{Generators.}
	\label{gen}
\end{figure}

There are many ways to view these models. Reshetikhin and Turaev \cite{RT} emphasized a categorical approach where the cups, caps and crossings of Figure~\ref{gen} are generators of a braided tensor category
whose morphisms correspond to {\it Morse diagrams} as in Figure~\ref{KnotTens} where the vertical direction on the page is the morphism direction and the diagram is cut up into cups, caps and crossings by horizontal lines.
In this approach the knot diagrams, arranged transverse to the vertical direction, are morphisms from a basic object $k$ to itself. The objects in the category are generated by $k$ and another object $V,$ but the main structure is in the morphisms. A diagram with ends, such as the crossing, is a morphism from $V \otimes V$ to itself. A functor from such a category of knots and tangles (diagrams with free ends) to modules over a ring $R$ will take the ``vacuum object" $k$ to $R$ and the ``free ends" object $V$ to a module $M$ over $R.$ The partition function or tensor contraction then appears as the result of the compositions of the linear maps to which the cups, caps and crossings go under the functor. In this way, the categorical formulation, the partition functions of statistical mechanics, the quantum amplitudes of quantum field theory and the tensor networks are all
unified in the categorical structure of these compositions. The reader can consult \cite{KD} for more details on this formulation of the quantum link invariants.

\begin{figure}[t]
	\centering
	\includegraphics[width=.6\textwidth]{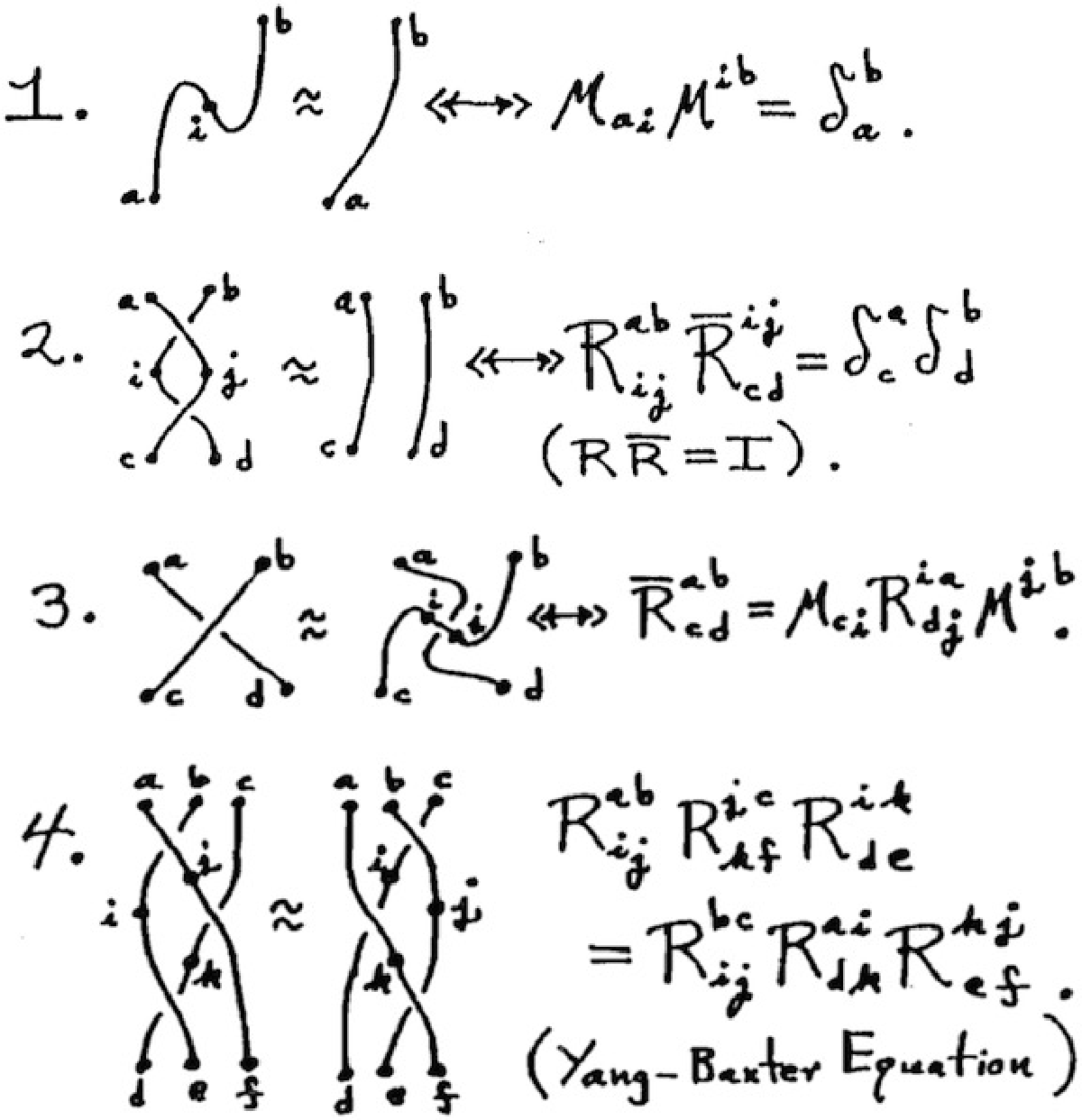}
	\caption{Tensor network equations for knot theory.}
	\label{tensor}
\end{figure}

\begin{figure}[t]
	\centering
	\includegraphics[width=.4\textwidth]{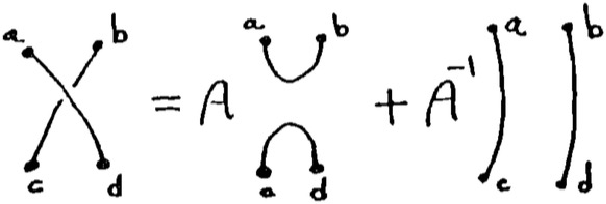}
	\caption{Tensor bracket expansion.}
	\label{tensorbracket}
\end{figure}

\begin{figure}[t]
	\centering
	\includegraphics[width=.4\textwidth]{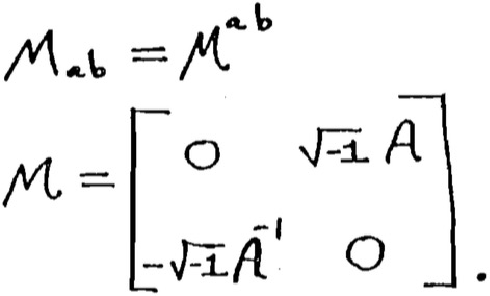}
	\caption{Cup-cap matrix.}
	\label{cupcap}
\end{figure}

\begin{figure}[t]
	\centering
	\includegraphics[width=.5\textwidth]{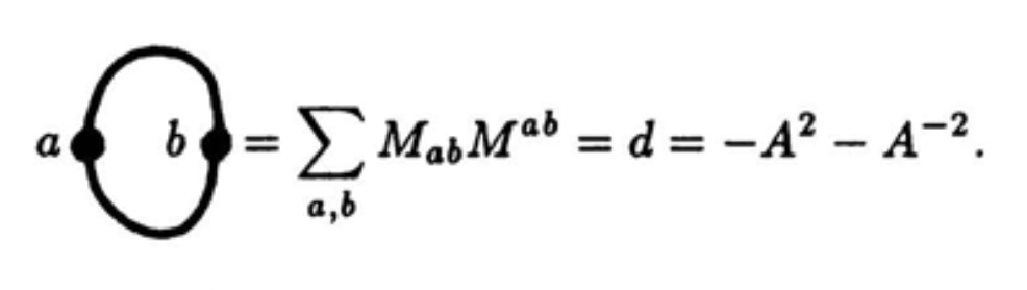}
	\caption{Loop evaluation.}
	\label{loop}
\end{figure}

In Figure~\ref{tensor} we show all the types of diagrammatic tensor equation that are needed to have an invariant of regular isotopy (second and third Reidemeister moves). These are generalizations of the 
Reidmeister moves to handle knots, links and tangles that are presented in the form of Morse diagrams. As the reader can see, we need (1) a cancellation of minima and maxima in unknotted lines, (2) a vertical version of the second Reidemeister move, (3) moves to shift an arc across a maximum or a minimum and the capacity to twist or shift an arc across a maximum or a minimum and (4) the Yang-Baxter equation. With matrix choices that satisfy these requirements for the cups, caps and crossings we can make knot invariants. For example, in Figure~\ref{tensorbracket}, we show how to translate the bracket expansion into a tensor equation.
This suggests that an appropriate choice of cups and caps will  provide the crossing tensor for the bracket. Indeed it is so and in Figure~\ref{cupcap} we show a $2 \times 2$ matrix $M$ such that $M^2$ is equal to the 
identity matrix. This $M$ is suitable for both the cup and the cap and gives the correct loop value for the bracket, as is shown in Figure~\ref{loop}. In this way it is not hard to use a small amount of algebra to produce 
a tensor network model for the bracket polynomial and hence for the Jones polynomial. The $R$-matrix made from this procedure is significant and occurs in the fundamental representation of the $SL(2)_{q}$ quantum
group. It is of great interest to examine knot theoretic invariants in the context of these quantum link invariants and tensor network models.

It should be remarked that what we have described for quantum link invariants, in many forms generalizes to the categorical formulations of topological quantum field theory, as defined by Michael Atiyah and Edward Witten. See \cite{Atiyah} for the first steps in 
this direction. A key point in topological quantum field theory is that it is the spaces and cobordisms of spaces that are objects and  morphisms in the associated categories. We will return to this point after discussing
the work of Edward Witten.

\section{Witten's Work, Quantum Field Theory and Vassiliev Invariants}
In 1988 Edward Witten discovered a quantum field theoretic approach to the Jones polynomial and its related invariants.
In \cite{Witten} Witten proposed a formulation of a class of 3-manifold
invariants and associated invariants of links in 3-manifolds via quantum field theory. 
He used generalized Feynman integrals in the the form  $Z(M,K)$  where
$$Z(M,K) = \int dAe^{(ik/4\pi)S(M,A)} W_{K}(A).$$
Here  $M$ denotes a 3-manifold without boundary and $A$ is a gauge field  (also
called a gauge potential or gauge connection)  defined on $M$.  The gauge field
is a one-form with values in a representation of
the  Lie algebra of $G$ for a specified Lie group $G.$ The group $G$ corresponding to this Lie algebra is said to
be the gauge group. In this integral the ``action"   $S(M,A)$  is taken to be the
integral over $M$ of the trace of the Chern-Simons three-form    
$$CS = A \wedge dA +(2/3)A \wedge A  \wedge A.$$  (The product is the wedge product of differential forms.) The term measuring the knot or link is 
$W_{K}(A),$ the trace of the holonomy of the gauge connection along the knot (product of such traces for links). The $k$ in the integral is an integral coupling constant.

$Z(M,K)$  integrates over all gauge fields modulo gauge equivalence (see
\cite{Atiyah} for a discussion of the definition and meaning of gauge
equivalence).

Witten's functional integral model of link invariants places them in the context of quantum field theory and quantum statistical mechanics.
In the form of this integral, this is the first time that we see the invariants expressed directly in terms of the embedding of the knot or link into three dimensional space.
All models described up to the point of Witten's work used diagrammatic representations for the topology. Witten's approach was a breakthrough into three dimensions and
into new relationships between topology and quantum field theory.
The formalism  and   internal logic of Witten's integral supports  the existence
of a large class of topological invariants of 3-manifolds and  associated
invariants of knots and links in these manifolds.

The three manifold invariants associated with this integral have been given rigorous
descriptions through the work of Reshetikhin and Turaev, Kirby and Melvin and Dror Bar-Natan
\cite{RT}, \cite{Kirby-Melvin}, \cite{Lickorish}, \cite{TL}, \cite{BarNatan}. The upshot of these descriptions is that the three-dimensional character of the invariants can be seen
via differential geometric expressions that arise in the perturbation expansion of the functional integral \cite{BarNatan}, but the original three dimensional vision of 
the integral remains problematic. 
Questions and conjectures arising from the functional integral
formulation are still outstanding. 
Specific conjectures about this integral take the form of just how it involves
invariants of links and 3-manifolds, and how these invariants behave in certain
limits of the coupling constant $k$ in the integral. Many conjectures of this
sort can be verified through the combinatorial and algebraic models. Some of the most perspicuous of these models use the work of 
Drinfeld \cite{Drinfeld,LOMI} on Hopf algebras to capture just the right context for the Yang-Baxter equation to reappear in the models in relation to the 
structure of the gauge groups. Drinfeld showed how solutions to the Yang-Baxter equation appear naturally in new algebras (the Drinfeld Double Construction)
that are directly related to the classical Lie algebras and to Hopf algebras more generally. In this way these algebraic results are deeply connected with the
quantum field theory.

The Witten integral can be explored via its perturbative expansion, just as is done in quantum field theory. This leads to relationships of the invariants defined by Victor Vassiliev \cite{Vassiliev}
with the coefficients in the perturbative expansion \cite{BarNatan,AF} and rapid development of Vassiliev invariants of finite type from this point of view \cite{CS}. Some of this development includes well-defined integral expressions 
for the Vassiliev invariants that go all the way back to the ideas of Gauss that defined integrals for linking numbers of curves in three dimensional space. In this way, the Witten integral did lead to a realization of the dream of a definition of the Jones polynomial in terms of an embedding of the knot or link in three dimensional space (instead of the combinatorial topology of the diagrams). The Vassiliev invariants also made clear,
via the work of Bar-Natan, Birman and Lin, how Lie algebras and their generalizations are fundamentally related to knot invariants. Up to the point of the introduction of the Vassiliev invariants there were two ways that 
Lie algebras entered the picture. Deformed  classical Lie algebras (aka quantum groups) figured in the work of Reshetikhin and Turaev to form knot invariants via categorical generalizations of state sums using
solutions to the Yang-Baxter equation. The deformations of the Lie algebras contained appropriate solutions to the Yang-Baxter equation. These techniques had turned out to be sufficient to reproduce on rigorous grounds the invariants that Witten defined by functional integration. But Lie algebras also figure in Witten's work via the choice of gauge group. Here it is a classical Lie algebra and a matrix representation of it that is chosen to produce a given invariant. The Vassiliev invariants give a unified point of view where the so-called weight systems for the Vassiliev invariant are computed from the Lie algebra and constitute initial data for 
integrating the invariant. The same initial data can be seen in the solutions to the Yang-Baxter equation that emerge from the quantum groups. With perfect hindsight one can see how the footprint of a Lie algebra -- the Jacobi Identity -- is related to topological invariance, and so one can draw the relationship of knot invariants and Lie algebras in a direct way that does not, in its logic, require either the quantum groups or the 
functional integrals. This is another aspect of this mathematics that deserves further understanding \cite{KP}.

\subsection{Vassiliev invariants and the Jacobi Identity}
\begin{figure}[t]
	\centering
	\includegraphics[width=.7\textwidth]{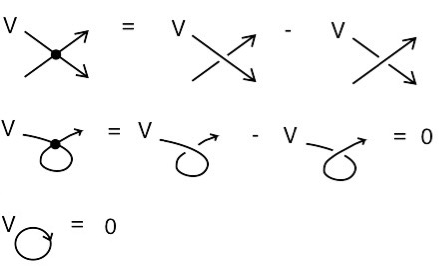}
	\caption{Vassiliev skein relation.}
	\label{vs}
\end{figure}

Link invariants are closely related with Lie algebras via the structure of solutions to the Yang-Baxter Equation that come from quantum groups (deformed Lie algebras) and from the gauge groups of the Chern-Simons-Witten theory. With this background it was eventually understood \cite{Birman and Lin,TS} how to relate Lie algbras directly to the knot theory via the Reidemeister moves. Here is a brief telling of that 
relationship. We shall say that $V$ is a  {\it Vassiliev invariant of finite type n} if $V$ satisfies the {\it Vassiliev skein relation} shown in Figure~\ref{vs} and  $V$ vanishes on all diagrams with more than $n$ nodes. The Vassiliev invariant is defined on knot and link diagrams that have the usual crossings but also have graphical nodes as illustrated in this figure. The skein relation says that the value of a diagram with a node is equal to the difference of the values of the two diagrams obtained by replacing the node by positive and by negative crossings. One can think of the diagram with the node as a kind of discrete derivative of the two diagrams with the crossings. In Vassiliev's viewpoint, the values of the graphical diagrams represent the differences between values assigned to different components of the space of all embeddings of knots.  It turns out \cite{Bar-Natan-Thesis,Bar-Natan} that when one makes the perturbative expansion of the Witten integral then finite type Vassiliev invariants appear as the coefficients of the inverse powers of the coupling constant. A similar result happens with the combinatorial version of the Jones polynomial if one makes a substitution of $e^{x}$ for the variable in the polynomial. Then the coefficients of $x^{n}$ are Vassiliev invariants of type $n.$\\

\begin{figure}[t]
	\centering
	\includegraphics[width=.7\textwidth]{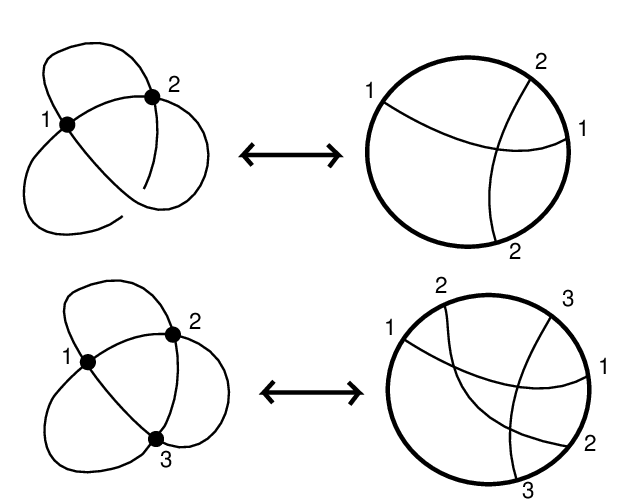}
	\caption{Chord diagram.}
	\label{chord}
\end{figure}

Experience mandates that one should look at these finite type invariants on their own grounds. Here is what happens: It follows from the difference equation of Figure~\ref{vs} that if $G$ represents a graph embedding  with $n$ nodes and $V$ has type $n,$  then $V_G$ is independent of the embedding of $G$ in three dimensional space. For diagrams this means that when $V$ has type $n$ and $G$ has $n$ nodes, then $V_G$ is independent of switching the crossings
in the diagram $G.$  For an example of this result, see Figure~\ref{chord} where we illustrate a diagram with two nodes. If we were computing a Vassiliev invariant of type 2, then the difference between the evaluation of the diagram shown in the figure and the one obtained by switching the crossing would be the value of the 3-noded diagram also shown in the figure. An invariant of type two will vanish on the 3-noded diagram. Hence the evaluation of the 2-noded embedding is independent of switching its crossings.
The evaluation of $V_G$ depends only on the graphical structure of $G$ defined by its nodes. It depends only on the structure of the chord diagram associated with $G$ that we define below.

\begin{figure}[t]
	\centering
	\includegraphics[width=.5\textwidth]{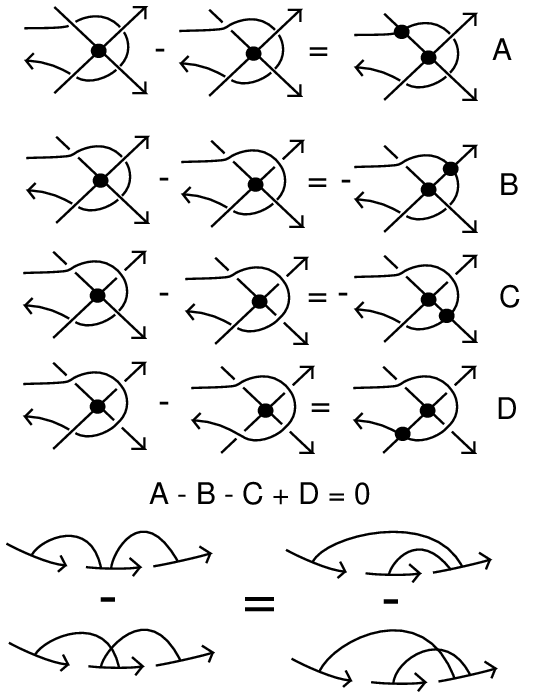}
	\caption{4 term relation.}
	\label{4t}
\end{figure}

\begin{figure}[t]
	\centering
	\includegraphics[width=.6\textwidth]{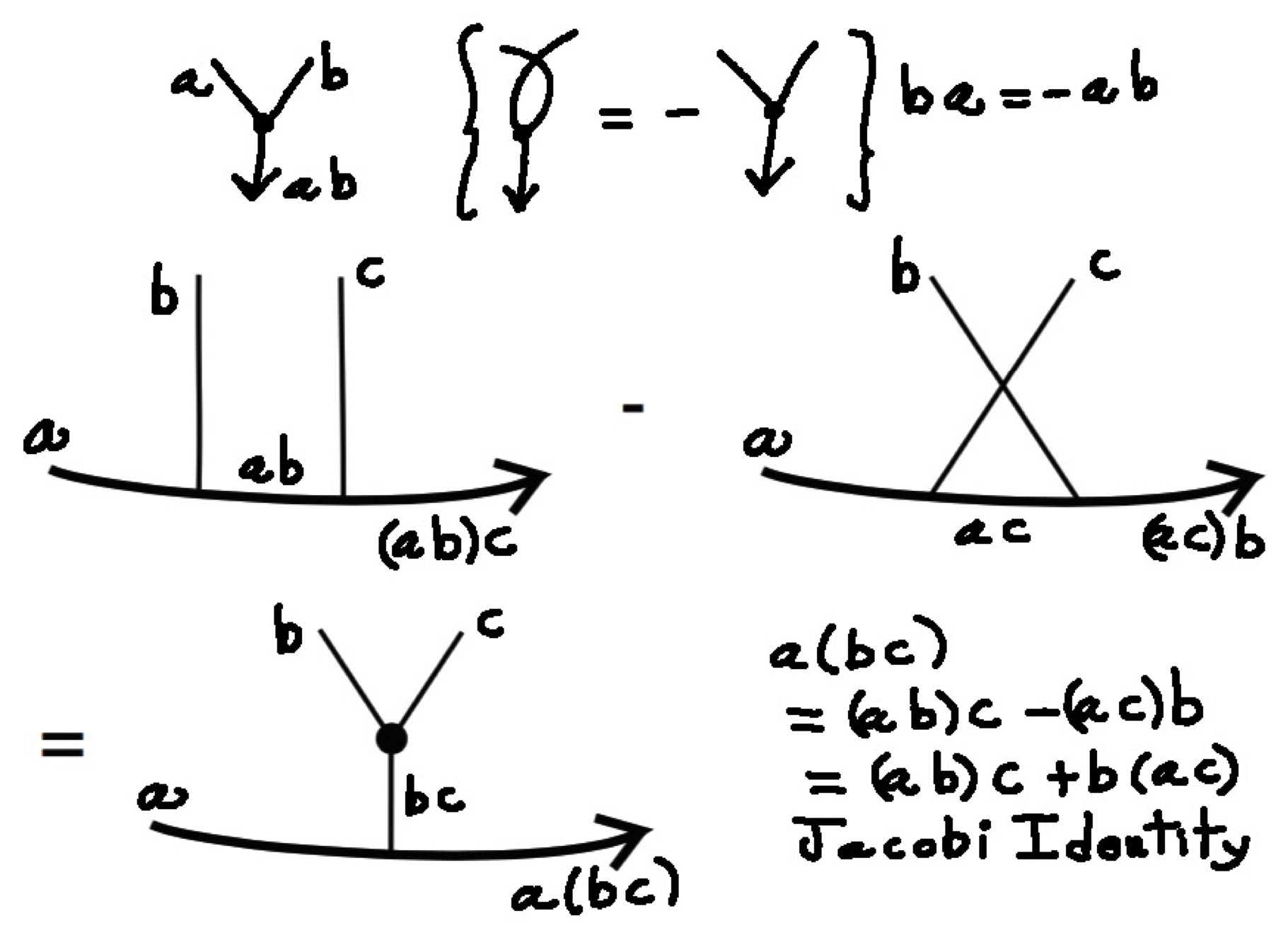}
	\caption{Lie algebra and Jacobi identity.}
	\label{jac1}
\end{figure}

\begin{figure}[t]
	\centering
	\includegraphics[width=.6\textwidth]{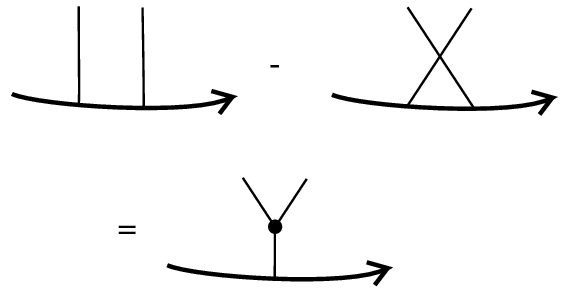}
	\caption{Jacobi identity.}
	\label{jac2}
\end{figure}

This evaluation is closely related to the Jacobi identity and to Lie algebras. One way to see this relationship is illustrated in Figures~\ref{chord},~\ref{4t},~\ref{jac1},~\ref{jac2}. In Figure~\ref{chord} we show how to encode the nodal information in a diagram $G$ in a so-called chord diagram. By taking a walk along $G$ one meets each node twice. This pattern of encounters is marked on a circle and chords are drawn between the pairs of appearances of the markers on the circle. In Figure~\ref{4t} we show how a relation on the chord diagram evaluations is obtained from the demand for invariance under the Reidemeister moves, coupled with the use of the Vassiliev skein relation. The reader will see four equations in this figure. Each equation is an instance of the Vassiliev skein identity. We start with a node and an arc that circles the node from underneath its edges. This is shown at the top left of the figure. This encircling arc passes under four points near 
the node. The first equation is a switching equation for the first point. The second equation is a switching equation for the second point, after the first node crossing has been switched. Proceeding in this way clockwise 
around the node, we obtain four equations. The diagram at the end has the encircling arc moving around the node from above. But by the (generalized) Reidemeister moves for these graph embeddings there is an equivalence between this last diagram and the very first diagram. This means that the sum of all of the left hand sides of these equations vanishes, and we are left with the statement that the sum of their right hand
terms is equal to zero, when evaluating them as Vassiliev invariants. As the reader can see from the figure, this is a sum of evaluations of embeddings. But when the diagrams shown have $n$ nodes for an invariant of type $n,$ then the resulting equation becomes the chord diagram relation shown at the bottom of the figure. This relation is called the {\it four-term relation} and is fundamental for the construction of 
Vassiliev invariants.

Figure~\ref{jac1} illustrates the definition of Lie algebra \cite{LA} and a diagrammatic representation of this definition. A Lie algebra $\cal{A}$ has a non-associative product, here denoted $ab$ for elements $a$ and $b$ of $\cal{A}$ with the 
properties
\begin{enumerate}
	\item Anticommutativity.  $ab = - ba$ for any $a$ and $b$ in $\cal{A}.$ 
	\item Jacobi Identity. $a(bc) = (ab)c + b(ac)$ for any $a,b,c$ in  $\cal{A}.$ 
\end{enumerate}
The second identity is the {\it Jacobi Identity} and can be regarded as the footprint of Lie algebra structure.
Lie algebras are ubiquitous in mathematics. One class of examples is the ring $\cal{M}(\cal{R})$ of $n \times n$ matrices over a commutative ring $\cal{R}$ with the Lie product taken to be the commutator of 
matrices $[A,B] = AB - BA.$ The reader will enjoy checking that the Jacobi identity is here satisfied as: $$[A,[B,C]] = [[A,B],C] + [B, [A,C]].$$  We have also that $[A,B] = -[B,A].$ With this, the usual Jacobi identity is seen as the equivalent form below:
$$[A,[B,C]] + [C,[A,B]] + [B, [C,A]]= 0.$$ 

Note that the Jacobi identity says that
left operation by an element of the Lie algebra satisfies the Leibniz Rule for products. Thus left multiplication is a derivation on the algebra. In Figure~\ref{jac1} we give a diagrammatic representation of these 
properties. Multiplication is represented by a trivalent node with the two upper legs labeled with elements $a$ and $b$, and the lower leg labeled by the Lie product $ab.$ Anticommutativity is represented by 
the trivalent node with crossed upper legs receiving a negative sign. The Jacobi identity appears as shown in the figure and the reader can follow the multiplications. Note that in the figure we find the formula
$(ab)c - (ac)b = a(bc).$ Since $- (ac)b = b(ac),$ this is the same as the Jacobi Identity as we have written it above. The diagrammatic advantage for writing it as we did is that the first diagram is two parallel lines
incident to a horizontal line. The second diagram is obtained by crossing the two vertical lines and the third diagram is obtained by running to the two vertical lines into a node that connects to the horizontal line. This means that the identity can be used in graphical networks by making local replacements. The outer edges of each term in the formula are the same and the vertical parts can receive the same algebra labels. In Figure~\ref{jac2} we illustrate the Jacobi identity with unlabeled diagrams. This identity can be transferred to a category of graphs or networks along with the anticommutativity so that we have Lie algebra diagrams as an abstract version of Lie algebras.

\begin{figure}[t]
	\centering
	\includegraphics[width=.6\textwidth]{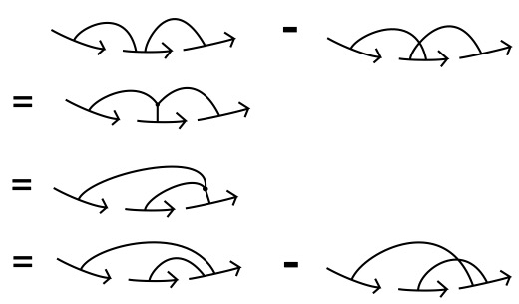}
	\caption{Diagrammatic proof of 4 term relation.}
	\label{diag0}
\end{figure}

In Figure~\ref{diag0} we show how Lie algebraic and chord diagrammatic structures come together. The figure is a {\it proof} showing that the 4-term relation we derived for Vassiliev invariants is a {\it consequence}
of the graphs being seen in a diagrammatic Lie algebra.  The difference in the left hand equation consists in two diagrams that differ only by a permutation. One has parallel lines.
The other has a crossing. This difference is replaced by a single network with a trivalent node. The trivalent node is moved by a planar isotopy to a new location and then opened up again by the reverse reading of the 
Jacobi identity. The result is the two terms of the 4-term relation on the right hand side of the equation. All of this can be clothed with specific algebra so that one obtains actual weight systems for Vassliev invariants 
from a multitude of choices of Lie algebras. This, in turn, gives rise to a host of invariants of knots and links (one must solve the problem of going from weight systems to actual invariants. See \cite{Vassiliev,Bar-Natan,Chmutov,Birman and Lin,GPV,KD,CS,TS}.). These 
invariants include the original Jones polynomial, the Homflypt and Kauffman polynomials, the quantum link invariants and more.

This relationship of the topology of knots and links via Reidemeister moves with the structure of Lie algebras is an extraordinary result. The full story of the relationship includes everything that we have indicated so far in this essay, from the origins of the Alexander polynomial and of the Jones polynomial, the quantum groups and statistical mechanics and the Witten functional integral. Yet, the story can be told in the very few lines by which we have described this relationship using the concept of 
Vassiliev invariants of finite type. The key role of the diagrammatic translations of language from knot and link diagrams, to chord diagrams, to network Lie algebra diagrams cannot be overemphasized.
With hindsight we can go from one diagrammatic language to the next and make a short line from knots and links to Lie algebras. The power of such translations of mathematical languages has not yet been fully tapped.

\section{Virtual Knot Theory and Extensions of the Jones Polynomial}
\begin{figure}[t]
	\centering
	\includegraphics[width=.5\textwidth]{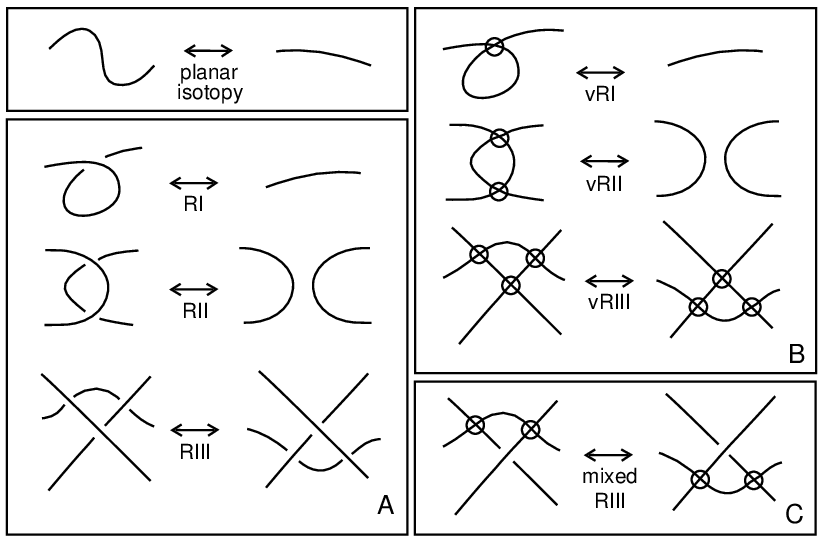}
	\caption{Virtual moves.}
	\label{vm}
\end{figure}

\begin{figure}[t]
	\centering
	\includegraphics[width=.5\textwidth]{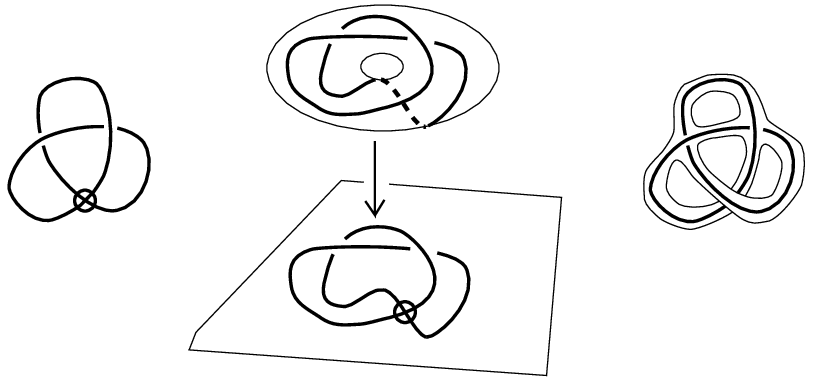}
	\caption{Surface representation of virtual knots and links.}
	\label{srep}
\end{figure}

In this section we describe a way to generalize the Jones Polynomial in a natural way to invariants of knots embedded in thickened surfaces, thus a way to 
have the Jones polynomial of a knot that is embedded in a three dimensional manifold other than the three sphere. The context of knots and links in thickened surfaces that we use
is called {\it virtual knot theory} \cite{VKT} and has a diagrammatic interpretation using an extra crossing type that I call a {\it virtual crossing}. In Figure~\ref{vm} we illustrate the virtual knot theory 
version of the Reidemeister moves. In Figure~\ref{srep} we illustrate how a knot diagram on a surface (representing a knot embedded in the thickened surface) corresponds to a planar diagram with 
virtual crossings. In that figure there is shown the projection of a toral knot diagram to the plane. Crossings that actually weave on the surface are projected to diagrammatic crossings
in the plane. However the projection also involves arcs that cross in the projection but are across the tube of the torus from one another. These are indicated in the planar diagram as
{\it virtual crossings}. The virtual crossings have no indication of over or under and are shown by a flat crossing of transversal arcs with a small circle drawn around the crossing point. A virtual diagram in the plane can be converted to a diagram on a surface in a standard way by forming a ribbon neighborhood of the diagram as shown in the right of Figure~\ref{srep}. Note that the ribbons at the virtual crossing are disjoint from one another. Once this standard neighborhood surface is formed, we add disks to its  boundary components and obtain a diagram of the corresponding knot or link on a surface.
If a knot is in a thickened surface and one can cut away handles on the surface in the complement of the knot, this is allowed in the virtual knot theory. In this way the virtual knot as embedding in a 
thickened surface is taken up to its embedding type in the surface and up to handle stabilization. Two virtual links are diagrammatically equivalent if and only if their surface representatives are stably equivalent in this sense. One can ask for the least genus surface in which a given virtual knot or link can be 
represented. Kuperberg \cite{Kuperberg} proved that there is a unique embedding type for each virtual link in its minimal genus surface. This gives a fully topological interpretation for this theory.
The moves in Figure~\ref{vm}  capture the stabilized equivalence relation and give us the freedom to represent a given virtual knot by a diagram with virtual crossings or to represent it on
a convenient surface. Two virtual links are said to be equivalent if one can be obtained from the other by the moves of Figure~\ref{vm}. In these moves the moves involving virtual crossings are designed
to insure that arcs of consecutive virtual crossings behave just as connections between their endpoints so that the arcs can be moved freely from place to place. 

We say that a virtual knot or link diagram is {\it classical} if, by virtual moves, it can be transformed to a diagram without virtual crossings, that being what we shall call a classical diagram.
It turns out that classical knot theory embeds in virtual knot theory in the sense that if $K$ and $K'$ are classical diagrams that are equivalent in the virtual context, then they are equivalent using only the 
classical Reidemeister moves \cite{Reidemeister}. In this way we see that the virtual knot theory is a genuine extension of the classical knot theory.

\begin{figure}[t]
	\centering
	\includegraphics[width=.4\textwidth]{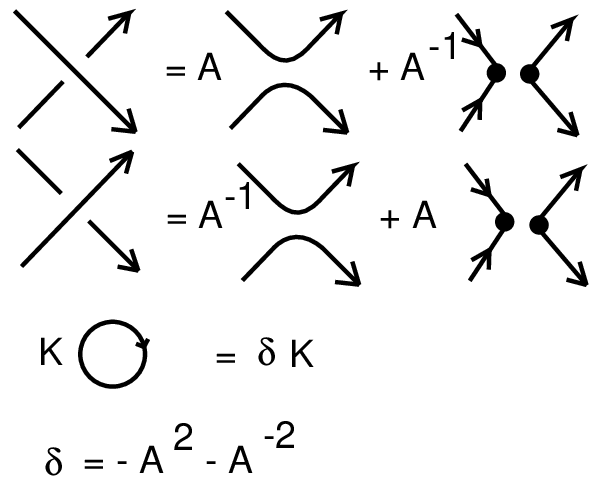}
	\caption{Arrow polynomial skein relations.}
	\label{arrow}
\end{figure}

\begin{figure}[t]
	\centering
	\includegraphics[width=.4\textwidth]{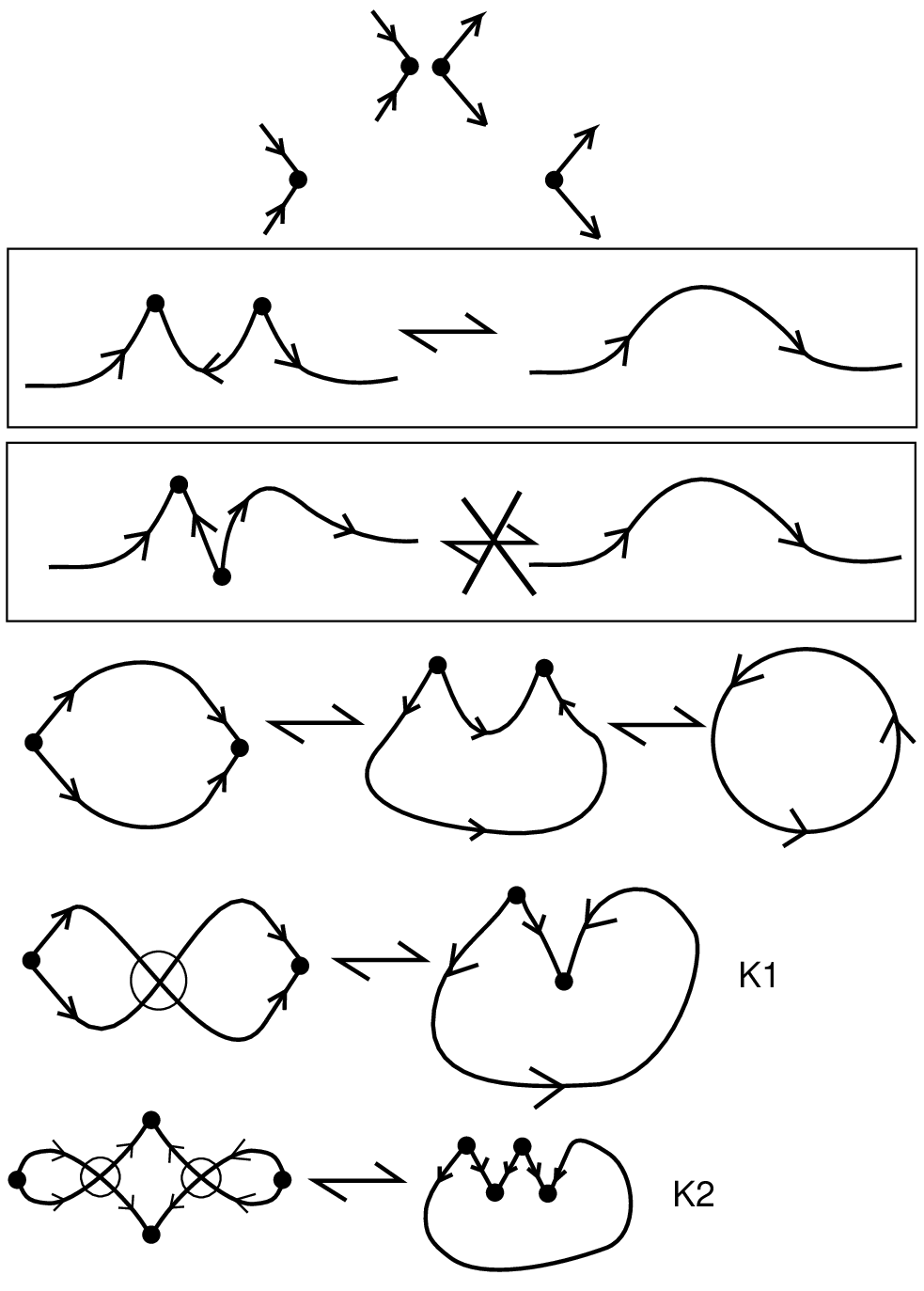}
	\caption{Arrow polynomial graphical reductions.}
	\label{arrowred}
\end{figure}

\begin{figure}[t]
	\centering
	\includegraphics[width=.4\textwidth]{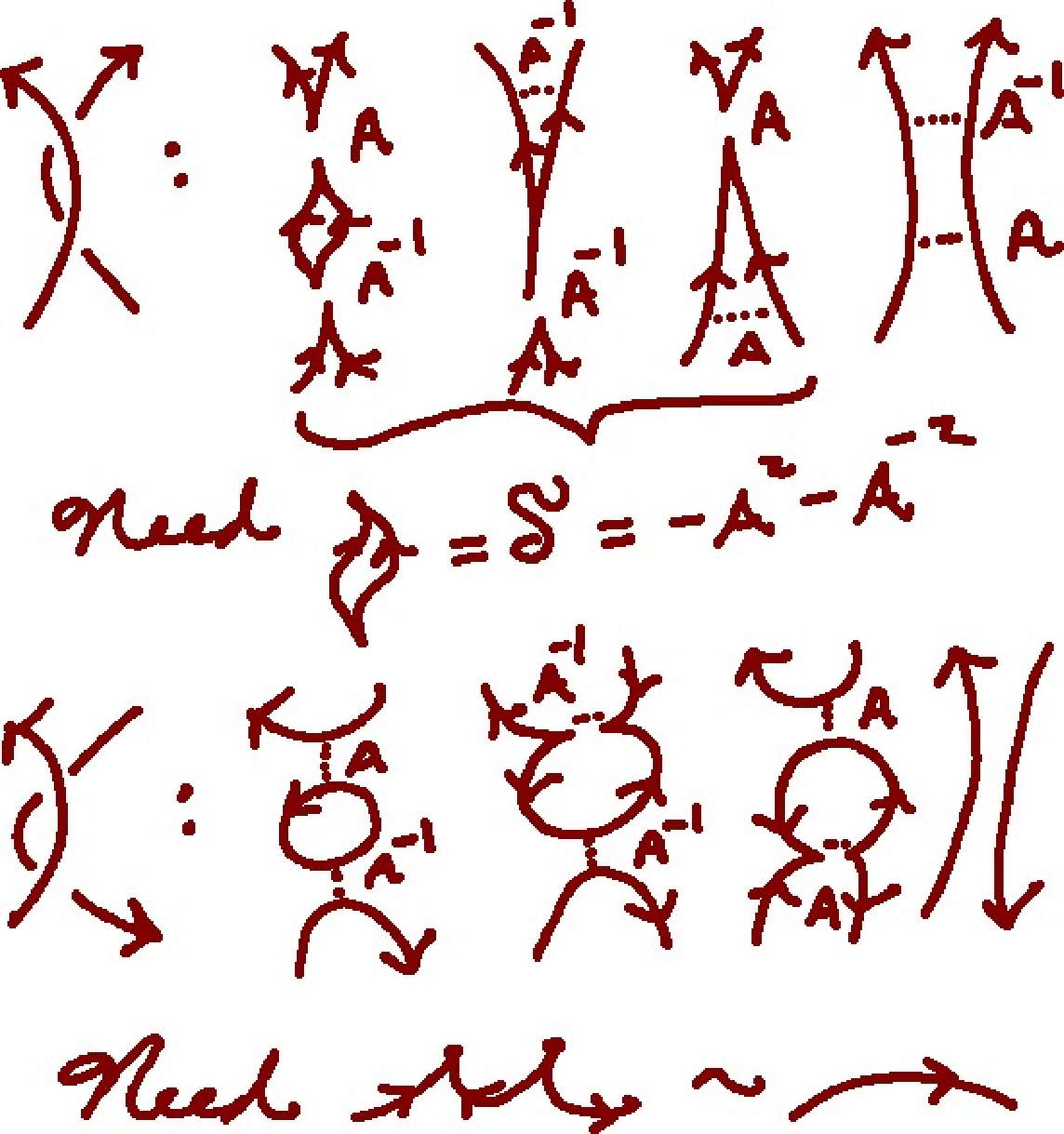}
	\caption{Arrow polynomial - Second Reidemeister move.}
	\label{arrowA}
\end{figure}

\begin{figure}[t]
	\centering
	\includegraphics[width=.4\textwidth]{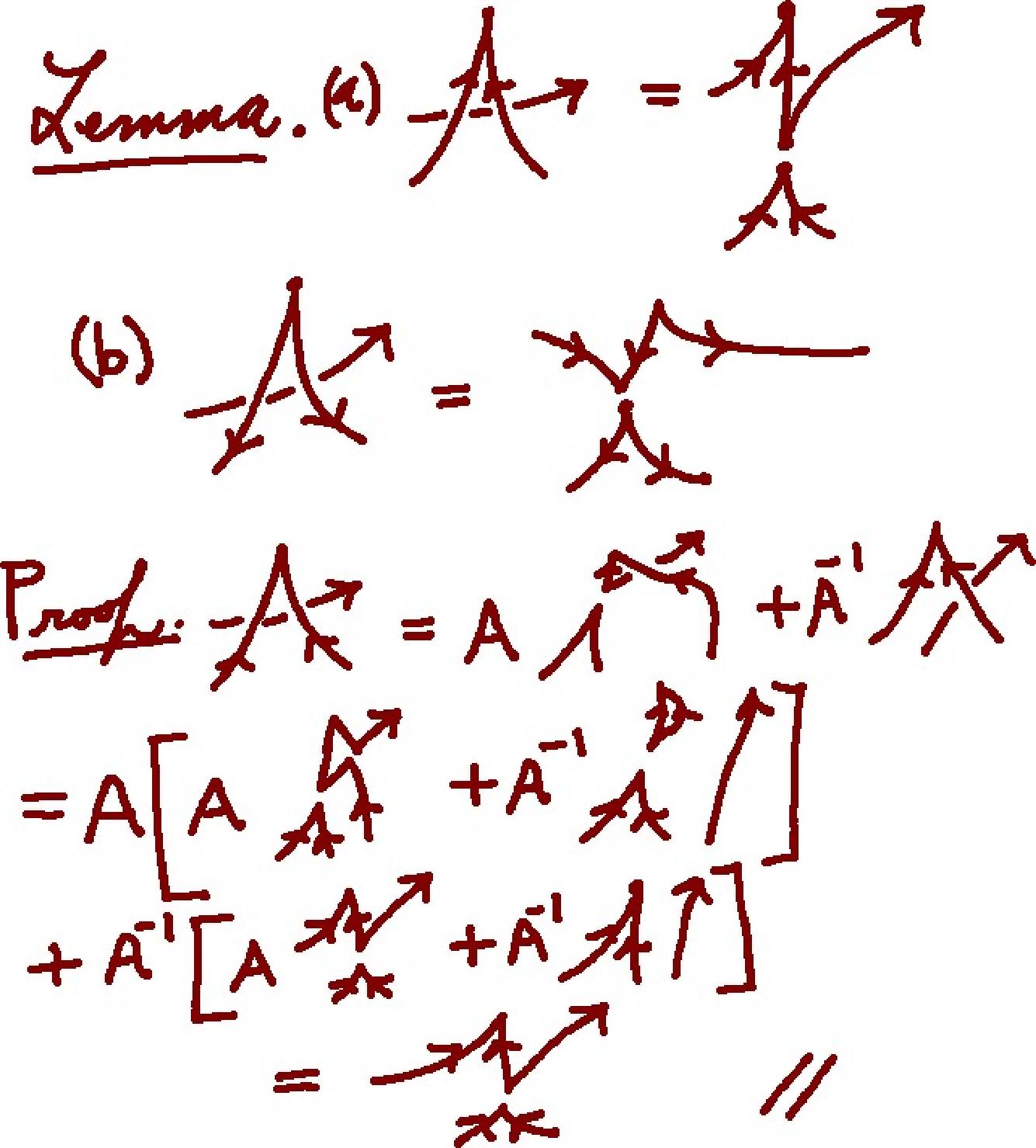}
	\caption{Arrow polynomial - The Zig-Zag Lemma.}
	\label{arrowB}
\end{figure}

\begin{figure}[t]
	\centering
	\includegraphics[width=.4\textwidth]{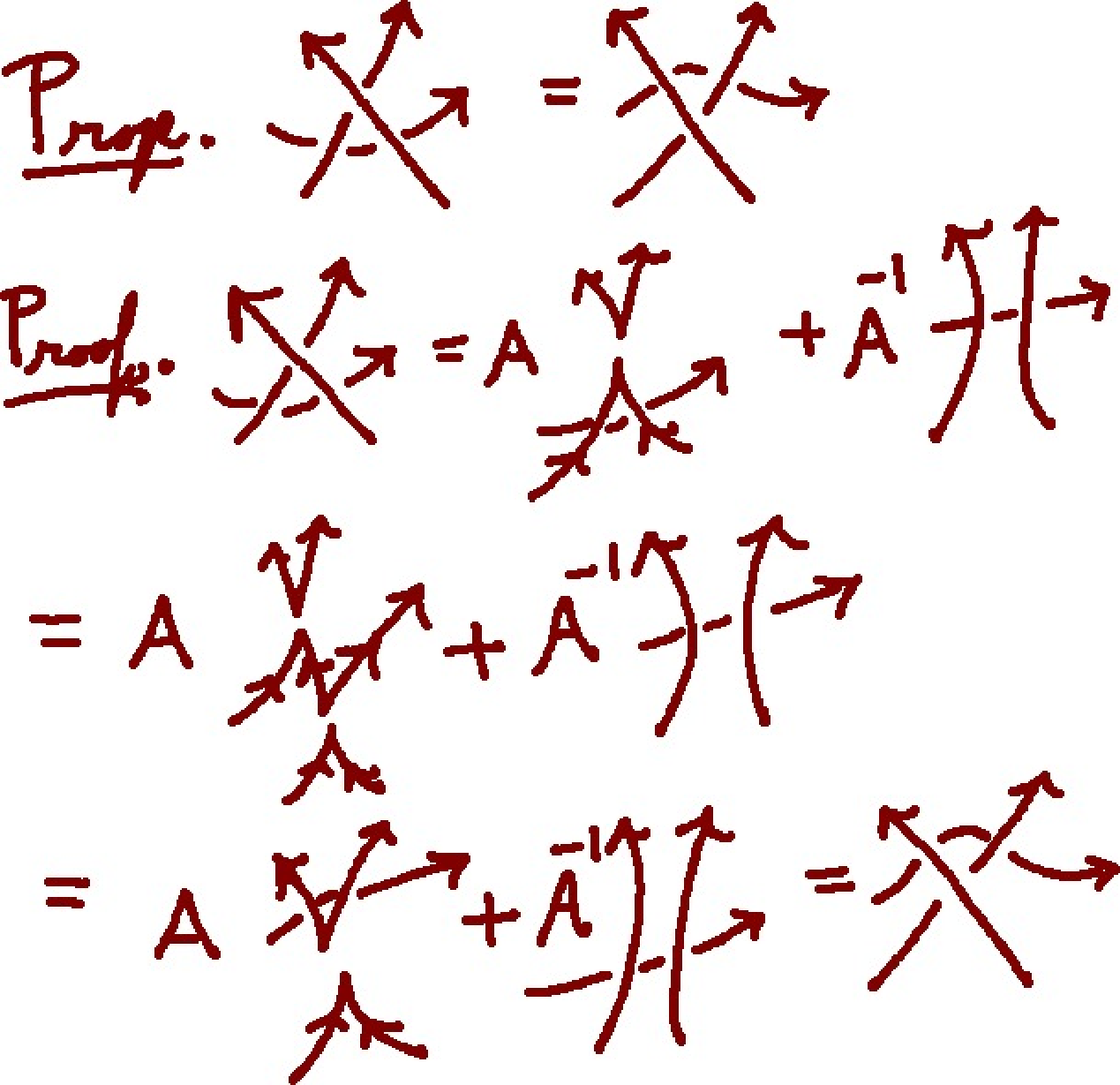}
	\caption{Arrow polynomial - Third Reidemeister move.}
	\label{arrowC}
\end{figure}

\begin{figure}[t]
	\centering
	\includegraphics[width=.4\textwidth]{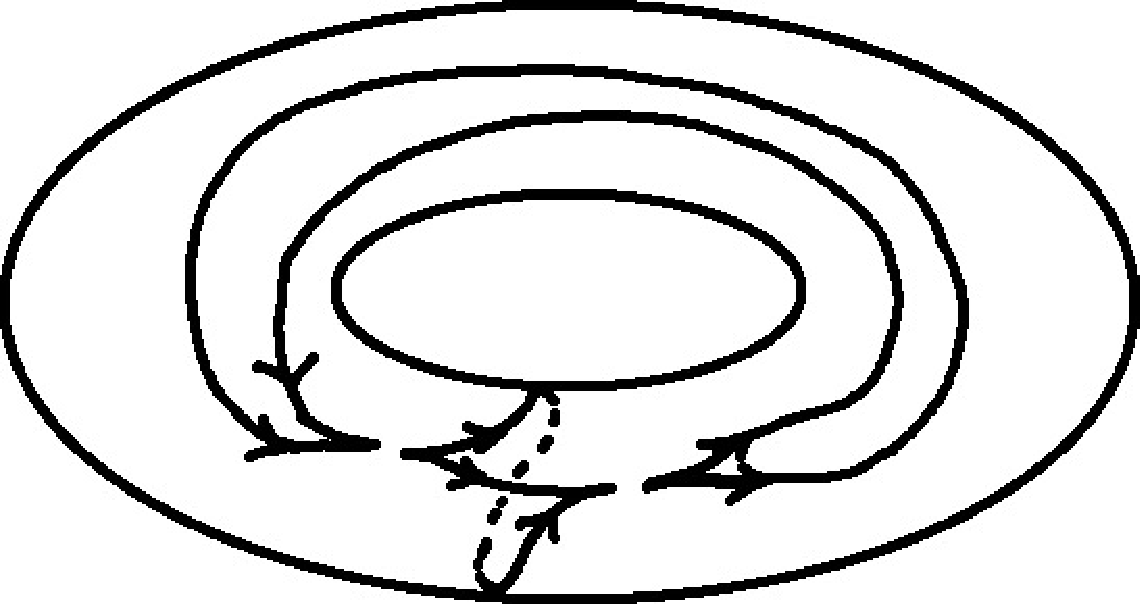}
	\caption{Zig-Zag state on a surface.}
	\label{zigzag}
\end{figure}

\begin{figure}[t]
	\centering
	\includegraphics[width=.4\textwidth]{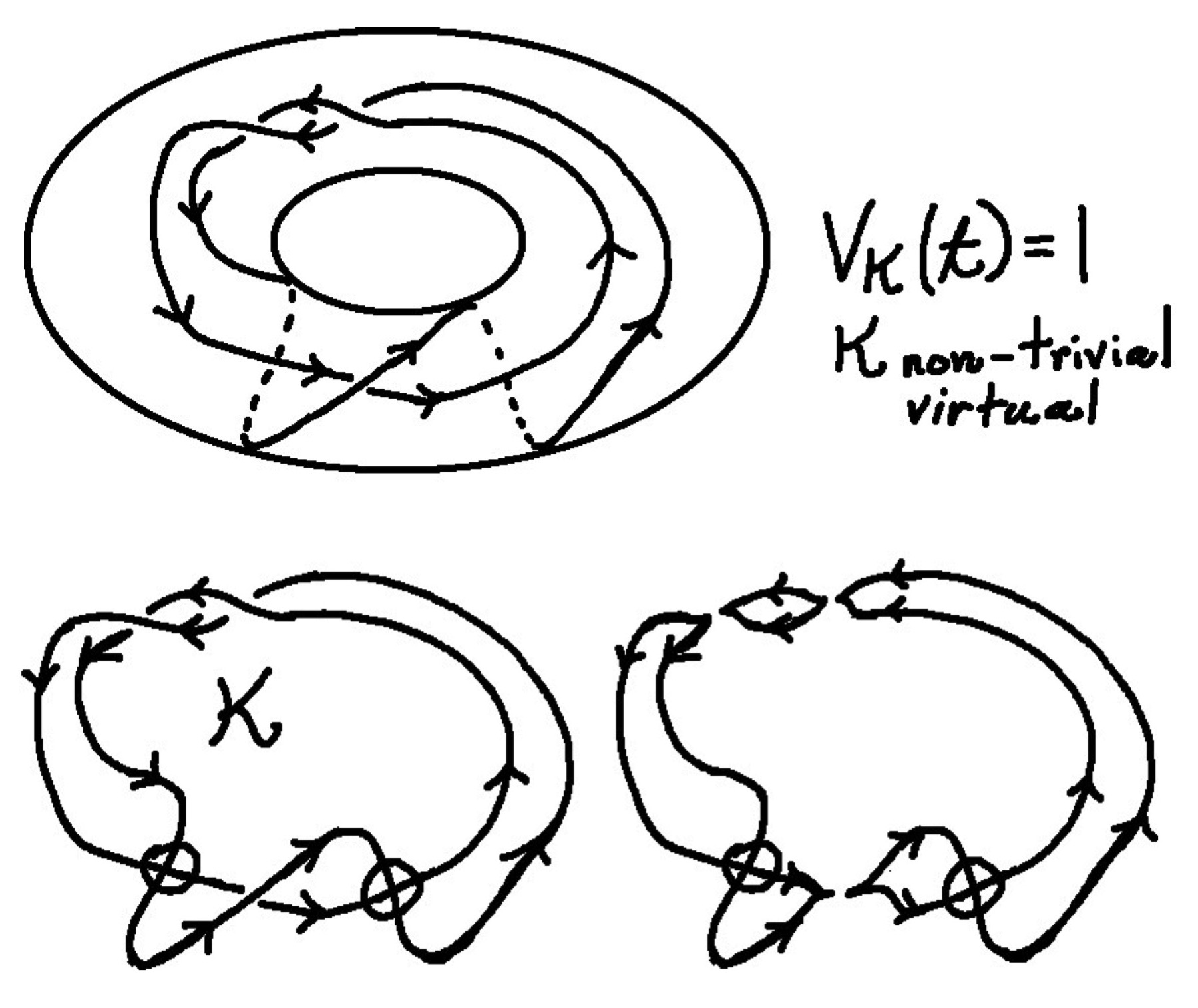}
	\caption{Zig-Zag state detects $K$ with unit Jones polynomial.}
	\label{zigzagdetect}
\end{figure}

The Jones polynomial in Kauffman bracket form generalizes directly to virtual knots by using the same expansion formula as the original bracket and evaluating loops in the same way 
without regard for the fact that the loops in a state expansion of the diagram may have virtual crossings. This natural extension of the Jones polynomial is of interest for classifying virtual knots, and it has
a special weakness. {\it There are non-trivial and non-classical virtual knots with unit Jones polynomial.} Such an example is the knot $K$ in Figure~\ref{zigzagdetect}. This weakness is of great interest
since it shows that by allowing the knot diagrams to leave the plane, the Jones polynomial definitely does not detect whether a knot is knotted. The corresponding problem for classical knots
remains an open problem at this writing. We do not know if there is a non-trivial classical knot diagram representing a non-trivial knot but having a unit Jones polynomial. It was this phenomenon of
non-trivial  virtual knots with unit Jones polynomial that led my original interest in formulating the virtual knot theory. There appears to be an analogy between the way non-planar graphs can be 
uncolorable (in the sense of the four color problem) while planar maps can (with a connectivity condition) always be colored and this phenomenon of knot detection for planar knots by the Jones polynomial
and the lack of detection when we leave the plane.  A better understanding of the plane is needed here.

The existence of knots that are not detected by this simple generalization of the Jones polynomial suggests that one should search for a better and more powerful generalization.
There is more than one way to find such invariants, but we shall describe here an invariant that we call the {\it Arrow polynomial} \cite{Arrow}. View Figure~\ref{arrow}. In this figure you see an oriented 
bracket expansion, and you see the strange way the orientation behaves when a crossing is smoothed in the ``B" way. In the figure for the positive crossing (first line of the figure) this smoothing has 
a coefficient of $A^{-1}.$ Call this the {\it disoriented smoothing}. As you can see the disoriented smoothing has two cusps with arrows going into one cusp and out of the other one.
It turns out that the resulting state sum will be a knot invariant just so long as we adopt the reduction rule shown in Figure~\ref{arrowred}. This reduction rule says that {\it two consecutive cusps on the same side of an arc can be canceled.}, as shown in the first box in the figure. The second box in the figure indicates that a zigzag of two consecutive cusps on two sides of the arc need not be canceled.
This means that reduced state loops with zigzags can be regarded as {\it new variables} for the Arrow polynomial. The argument that verifies this fact works for diagrams on surfaces, for virtual diagrams and for classical planar diagrams.

The diagrammatics for the proof of invariance of the Arrow polynomial under the Reidemeister two and Reidemeister three moves is given in Figures~\ref{arrowA}, \ref{arrowB}, \ref{arrowC}.
The reader should have no difficulty deducing invariance under the second Reidemeister move from Figure~\ref{arrowA}, using the reduction rule that allows cancellation of consecutive cusps on the same side of an arc.
In Figure~\ref{arrowB} is given the proof of the ``Zig-Zag Lemma" that shows how an undercrossing arc at a cusp can be pulled out from under the cusp, acquiring a zig-zag. This lemma is them applied to the skein
expansion of the third Reidemeister move to first produce such zig-zag and then remove it in the course of transforming the arc underneath a crossing as is desired. This proof is a generalization of the invariance 
proof for the original bracket polynomial.

Irreducible zigzag states do not occur for classical diagrams, but zigzags naturally occur when the knots are in 
thickened surfaces. Figure~\ref{zigzag} shows how a zigzag can occur in a state on a torus. Figure~\ref{zigzagdetect} indicates a zigzag state for a knot $K$ with unit Jones polynomial
The Arrow polynomial detects the non-triviality of this $K$.

There are many interesting properties of the Arrow polynomial, and we do not yet know just how strong an invariant it is. The Arrow polynomial 
is, in its way, a justification for the definitions of virtual knot theory. It is a generalization of the Jones polynomial that is natural when you start with oriented diagrams and use the bracket expansion. 
A person following the logic for constructing the Arrow polynomial would be led inevitably to invent virtual knot theory. Things do not always happen in such straight lines. A reader of my book
``Knots and Physics'' \cite{KP} will discover that I almost invented the Arrow polynomial in those days when I wrote the book (and actually when I was working out the bracket in Venice in the fall of
1985). But not quite, and so it was only later, after virtual knot theory had emerged,  that Heather Dye and I found the Arrow polynomial, and it was in fact independently found by Naoko Kamada and Miyazaki  with a 
slightly different viewpoint \cite{KM}. By this excursion into invariants of knots in three-manifolds we have found very beautiful generalizations of the Jones polynomial, and many new questions. 
It is to be hoped that there are more insights of this nature in the exploration of other three manifolds and their knot theories.

\section{Khovanov Homology}
We now discuss Khovanov's discovery \cite{Kho} of a  categorification of the Jones polynomial. 
We have already introduced the Khovanov Category  $Cat(K)$ of a knot or link diagram $K$ in Figure~\ref{cube} and the cube category in Figure~\ref{cubecat}.
Now we can face this question: {\it How can one extract topological information about the knot diagram $K$ from its category $Cat(K)$?}  Khovanov succeeded in doing this, and I want to introduce his 
construction and its relationship with the bracket polynomial and the Jones polynomial by exploring this question from our view of the category $Cat(K)$.
In this section we discuss a point of view about Khovanov Homology that was developed by Dror Bar-Natan. This section is a sketch of that point of view. The reader can find more detail in the papers
\cite{BN1,BN2,LKKho}.

\begin{figure}[t]
	\centering
	\includegraphics[width=.45\textwidth]{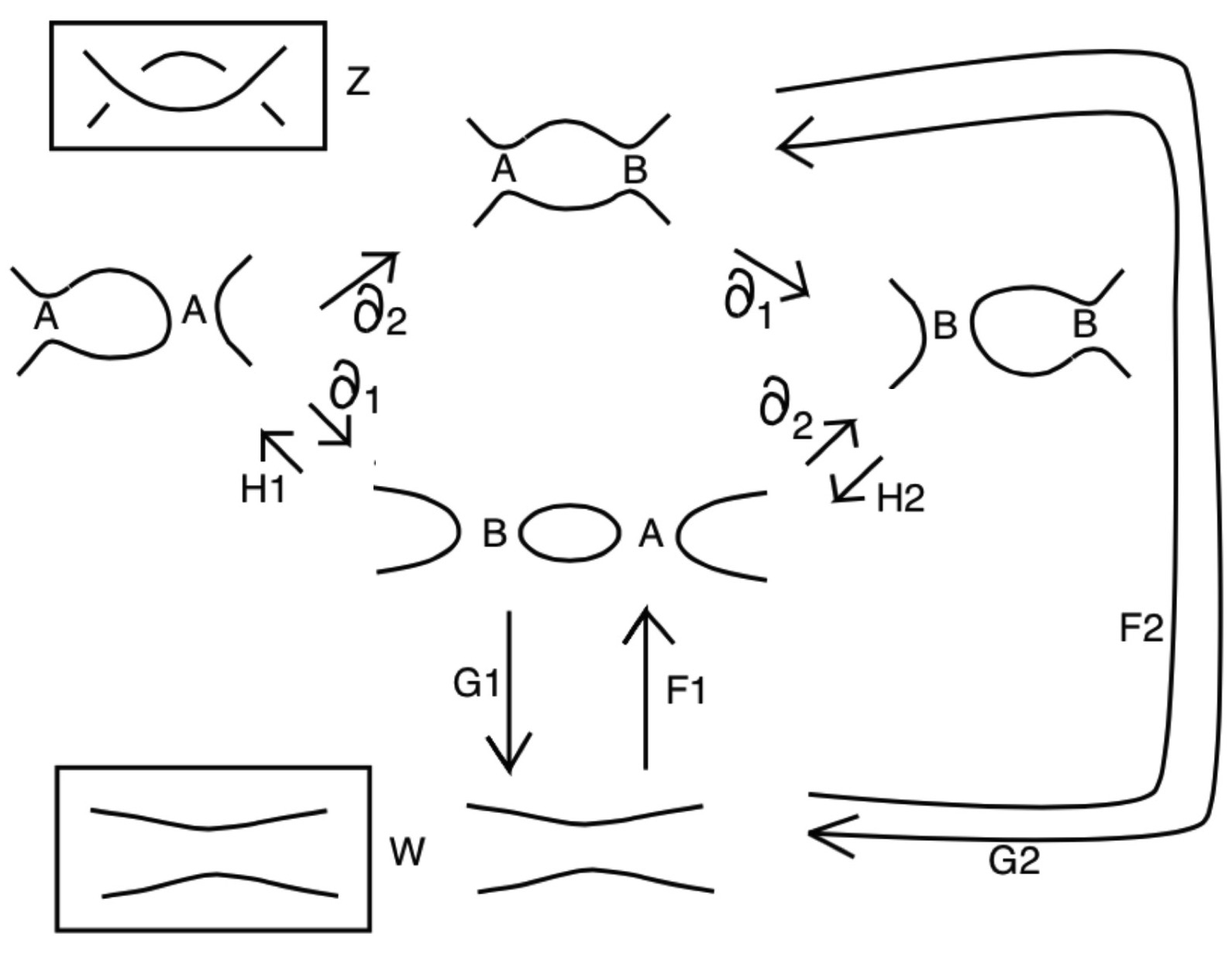}
	\caption{Khovanov category for Reidemeister two move.}
	\label{diag}
\end{figure}

Consider how $Cat(K)$ will change if we apply a Reidemeister move to $K$. In particular, consider the change that results from a type two Reidemeister move. We indicate this change in the category with 
the diagrams in Figure~\ref{diag}. In this figure the arrows from a state with an $A$-smoothing to a state with a $B$-smoothing are labeled $\partial_{1}$ and $\partial_{2}.$ The smoothings go in order from left to
right in the diagram and $\partial_{1}$  denotes operating on the left smoothing, while  $\partial_{2}$ denotes operating on the right smoothing. For the upper boxed diagram with two crossings, there are four local smoothings that are linked by these arrows. The lower boxed diagram is the result of a  type two Reidemeister move applied to the upper boxed diagram. Here there are no crossings, so the category of the lower boxed diagram can be represented locally by one already-smoothed diagram as is shown. In order to have a map between the lower category and the upper category we seem to need to map single objects (states) in the lower category to multiple objects in the upper category. Furthermore we would like to have a description of a map $F_{1}$ from the lower category diagram with parallel arcs to a diagram that is identical with it but has 
a circle in between the parallel arcs. A better way to think about the category is needed!

\begin{figure}[t]
	\centering
	\includegraphics[width=.45\textwidth]{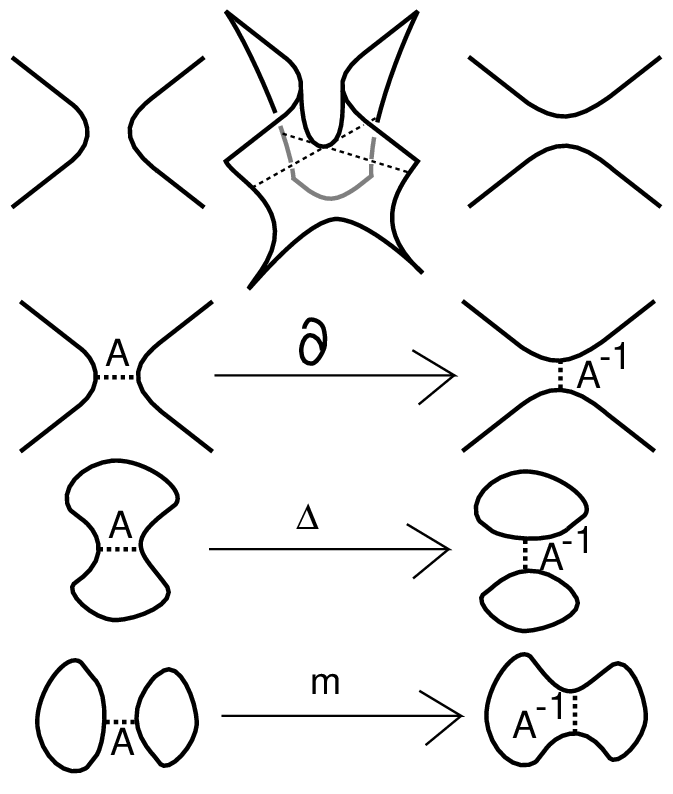}
	\caption{Saddle point resmoothing morphism.}
	\label{saddle1}
\end{figure}

\begin{figure}[t]
	\centering
	\includegraphics[width=.6\textwidth]{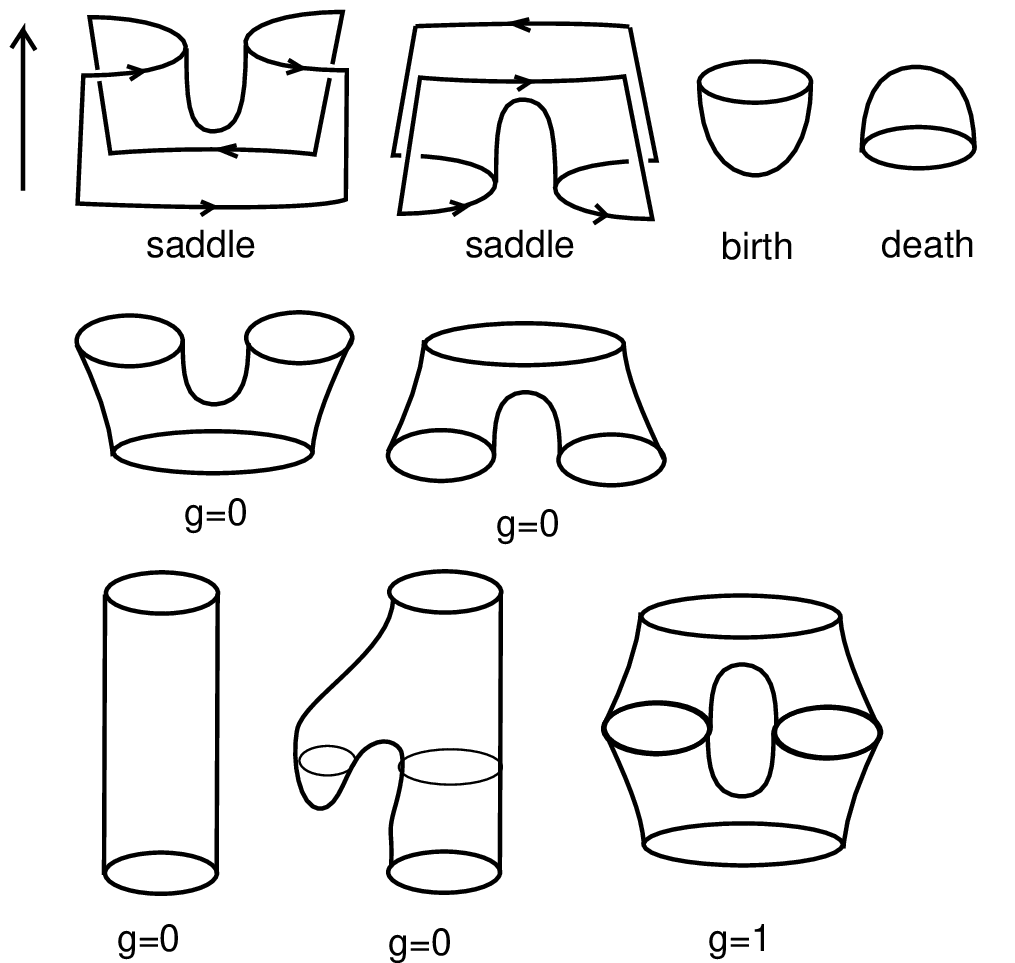}
	\caption{Saddles, births and deaths.}
	\label{saddle}
\end{figure}

Please view Figure~\ref{saddle1}, Figure~\ref{diag},  and Figure~\ref{saddle}. Here we indicate a way to think about the morphisms in the category. A resmoothing can take one loop to two loops, or two loops to one loop. We can think of these loops as the boundary ends of a surface that connects them by a saddle point as illustrated in the figures. From this point of view each morphism in the Khovanov category can be regarded as a surface that begins with one state as its ``left" boundary and has the target state as its ``right" boundary. In this way of thinking it is also possible for an individual circle to be created or destroyed by taking it through a minimum or a maximum in the course of the sections of the surface. With this way of looking at the category, we can construct the map $F_{1}$ in the previous paragraph, by allowing the ``birth" of a circle via a minimum so that the circle that is seen in the upper part of the diagram at the bottom is the boundary of a bowel and one sees the bowel as a creation process for the circle.

\begin{figure}[t]
	\centering
	\includegraphics[width=.7\textwidth]{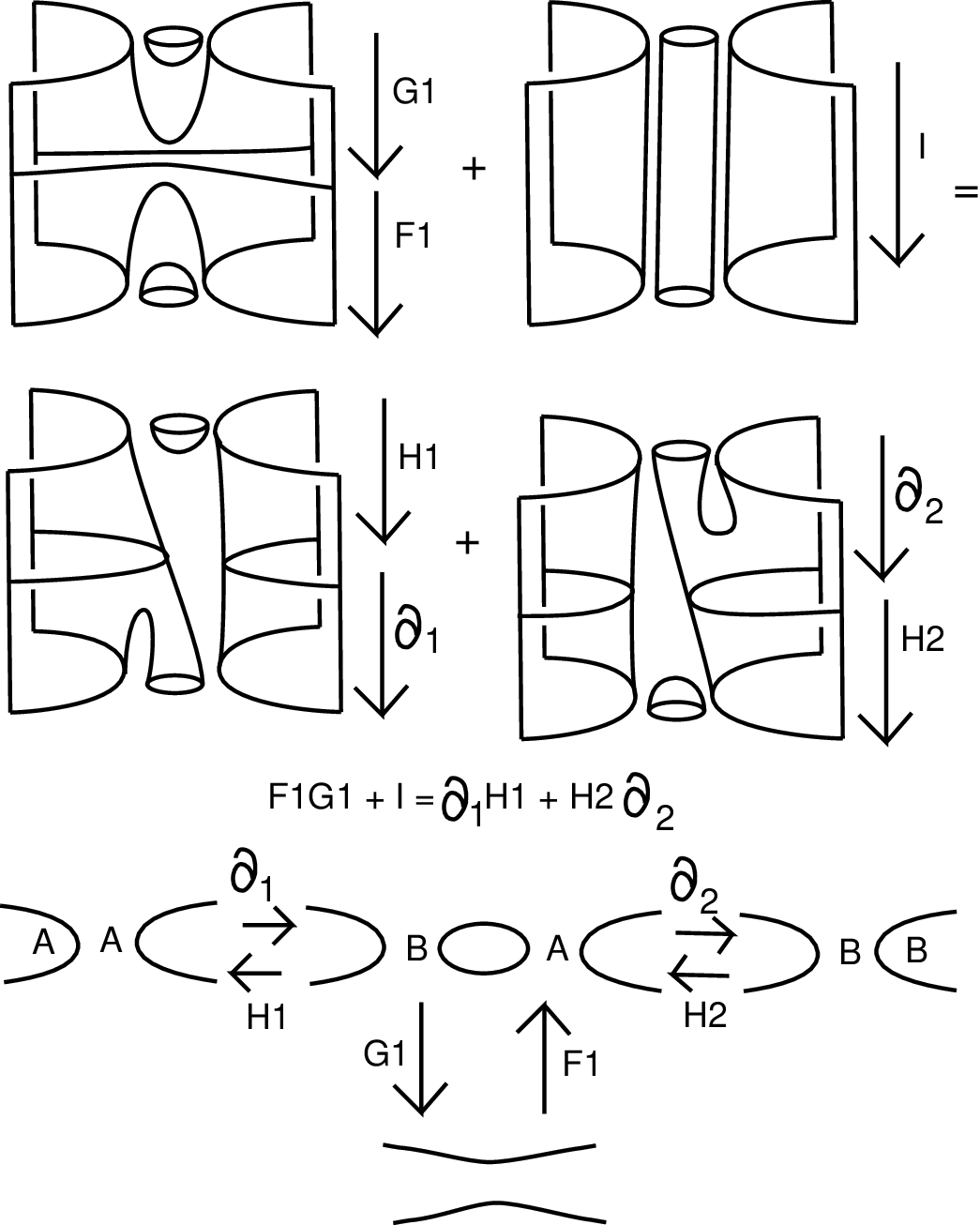}
	\caption{Homotopy for Reidemeister two move.}
	\label{diag1}
\end{figure}

\begin{figure}[t]
	\centering
	\includegraphics[width=.6\textwidth]{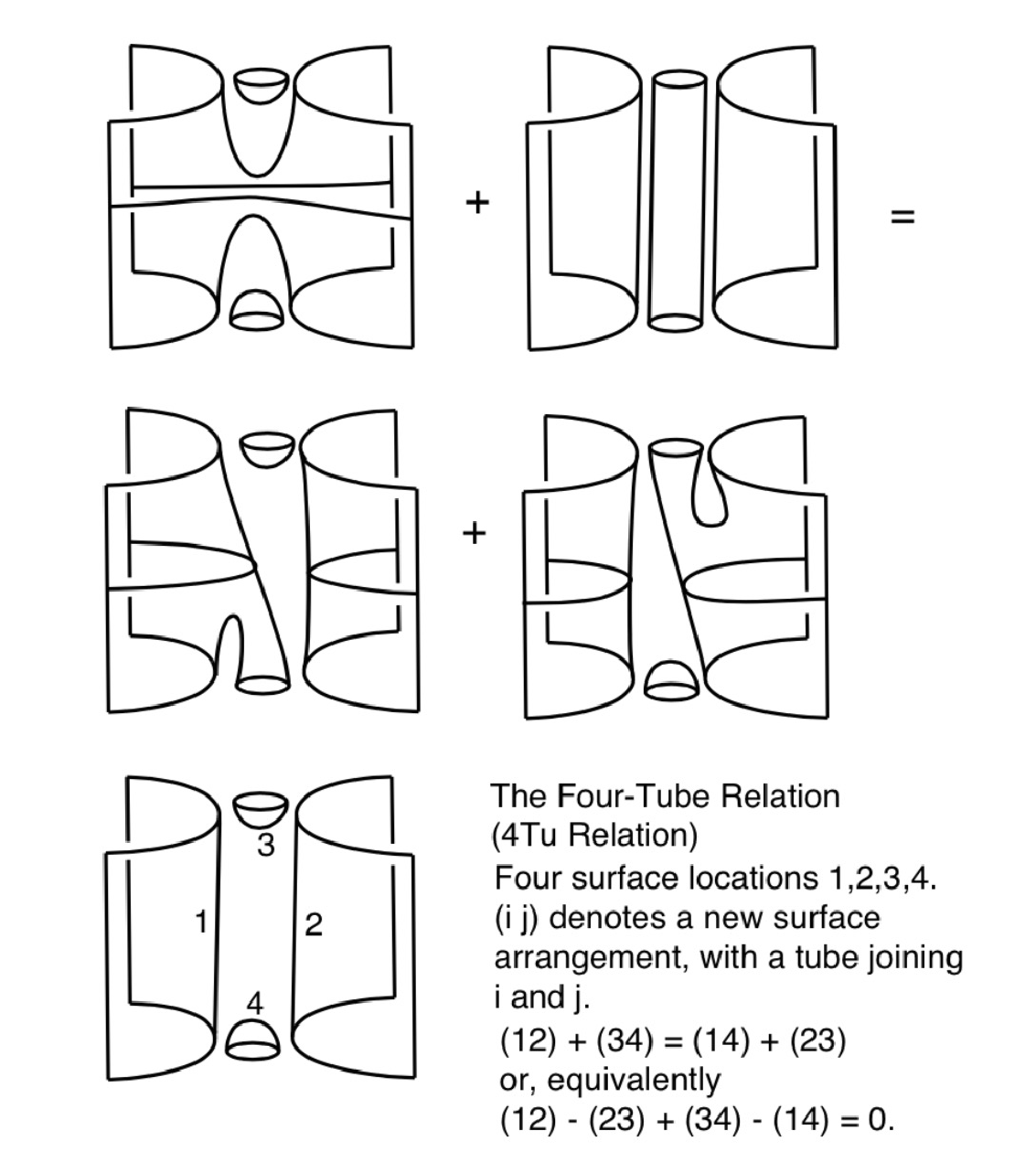}
	\caption{4 tube relation.}
	\label{diag2}
\end{figure}

\begin{figure}[t]
	\centering
	\includegraphics[width=.4\textwidth]{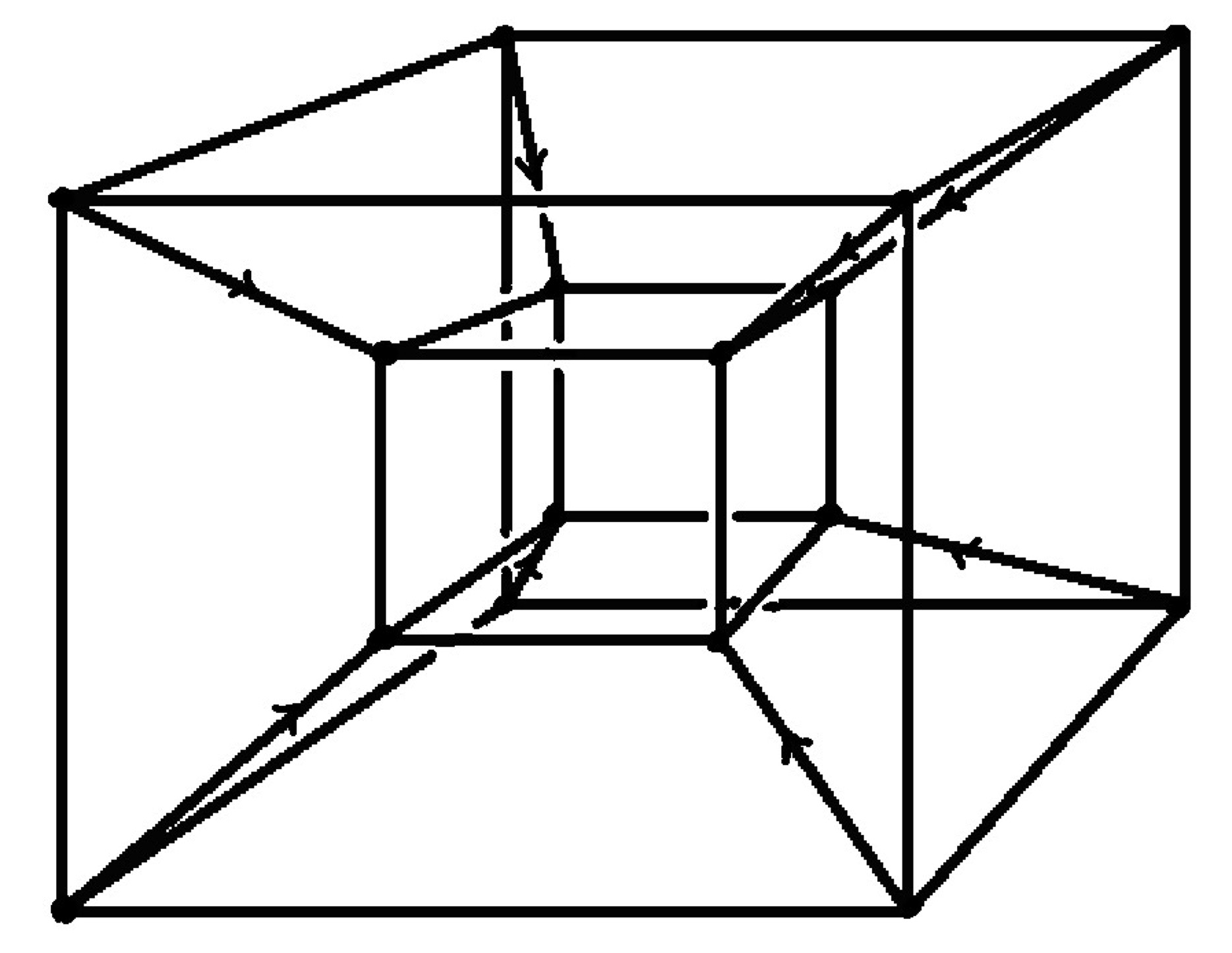}
	\caption{4-cube category.}
	\label{4cube}
\end{figure}

\begin{figure}[t]
	\centering
	\includegraphics[width=.4\textwidth]{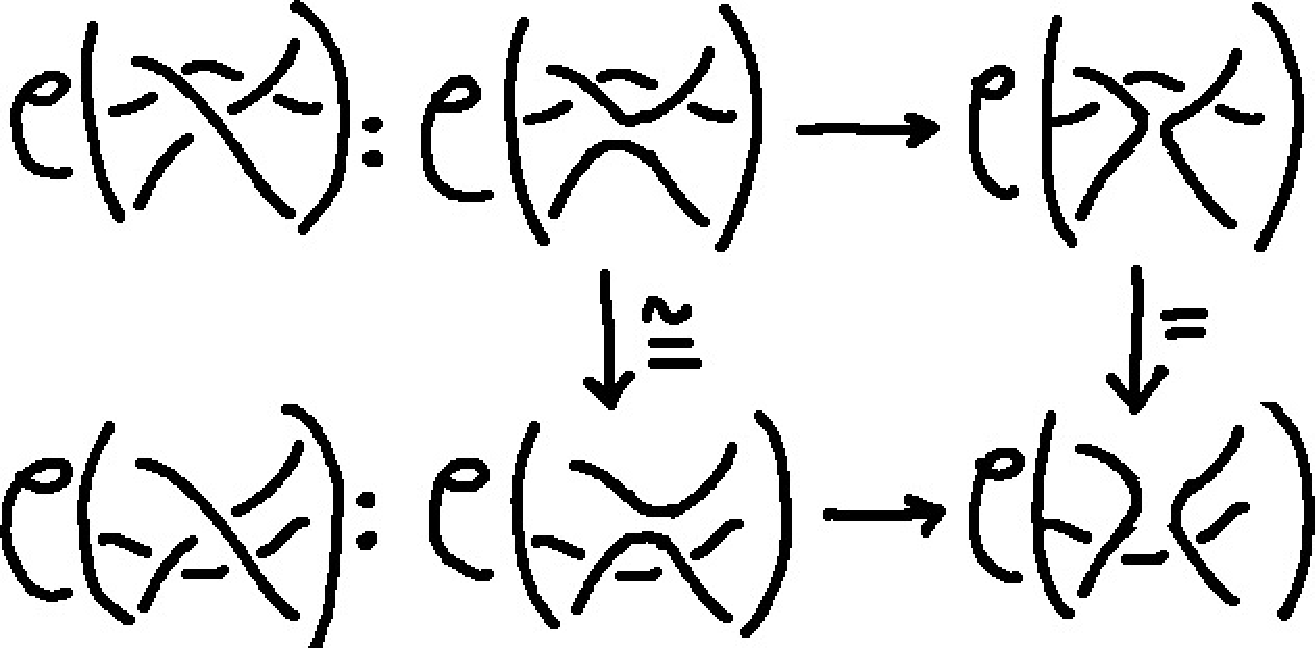}
	\caption{Categorical mapping cones.}
	\label{mcone}
\end{figure}

But there is still a problem in comparing the two categories. It is the problem that the lower category wants to map objects to multiple objects in the upper category. We can solve this problem by allowing multiple objects to become single objects. We are familiar with this idea in algebra where we take the direct sum of algebras. By the same token we would like to take the direct sum of the two objects that lie directly above the 
parallel arcs in Figure~\ref{diag}. This can be done by forming a new category that we shall call the  {\it Bar-Natan Cobordism Category} \cite{BN1,BN2} where all the states with the {\it same number} of $B$-smoothings are amalgamated into 
one direct sum object. Our original category now is rearranged and has the form $$C^{0} \longrightarrow C^{1} \longrightarrow C^{2} \longrightarrow \cdots$$ where $C^{k}$ is the direct sum object corresponding to all the states that have $k$ $B's$. Look again at Figure~\ref{cube} and imagine amalgamating the horizontal rows into a ``direct sum" of the states that are in the row. A sequence of objects and arrows as we have drawn above reminds us of a chain complex in algebraic topology. The Bar-Natan Cobordism Category is an abstract chain complex. The objects and morphisms are not maps of modules, but functors can take them to maps of modules and make honest chain complexes from them. We will continue to look at the Bar-Natan Cobordism Category in its simple abstract form. One can add maps abstractly and also subtract them. The sums of arrows from $C^{k}$ to $C^{k+1}$ can be seen to be good boundary maps in the sense that the compositions $C^{k} \longrightarrow C^{k+1} \longrightarrow C^{k+2}$ are zero with appropriate assignments of signs. One way to assign the signs is to order the crossings in the original knot or link diagram $K$ so that a state corresponds to an object in the cube category via a sequence of letters $A$ or $B$ corresponding to the local smoothings. For example in Figure~\ref{cubecat} we could have $[A,A,B]$ stand for one of the states of the trefoil knot in Figure~\ref{cube}. Then we would have a resmoothing map from $[A,A,B]$ to $[A,B,B]$ and we want to see in the Bar-Natan Cobordism Category whether to add it or subtract it. One answer that works is {\it if you are resmoothing an $A,$ assign $(-1)^{t}$ where $t$ is the number of $A's$ that precede the $A$ you are smoothing in the given state.} Thus for 
$[A,A,B] \longrightarrow [A,B,B]$ we would use a minus sign. We let $\partial : C^{k} \longrightarrow C^{k+1}$ be the sum of all the maps in $Cat(K)$ from states with $(k) B's$ to $(k+1) B's$ with these signs.  Then 
$\partial \partial = 0$ in the Bar-Natan Cobordism Category,  and we have an abstract chain complex.

Two chain complexes can be compared by {\it chain homotopies}. Two chain maps $f,g: C^{*} \longrightarrow D^{*}$ are said to be chain homotopic if there is a mapping $H:C^{*} \longrightarrow C^{*-1}$ such that 
$\partial H + H \partial = f + g.$ We will not worry about signs here. Note that $\partial : C^{*} \longrightarrow C^{*+1}$ is the boundary map for the chain complex so that $\partial \partial = 0.$ Here we use upper indexing as in the previous paragraph. Two chain complexes $C^{*}$ and $D^{*}$ are {\it chain homotopy equivalent} if there are chain maps  $f:  C^{*} \longrightarrow D^{*} $
and $g:D^{*} \longrightarrow C^{*} $ so that each composition $fg$ and $gf$ is chain homotopic to the identity. 
Complexes that are chain homotopy equivalent have the same homology groups where the reader will recall that the $k$-th homology group
$H^{k}(C^{*})$ is the kernel of the boundary map from $C^{*}$ modulo the image of the boundary map from $C^{*-1}$.

All these concepts about chain homotopy go over to the Bar-Natan Cobordism Category (except for the calculation of kernels and images).  Thus we can consider the chain homotopy class of the Bar-Natan Cobordism Category of 
a knot or link $K$. Let $KhoCob(K)$ denote the Bar-Natan Cobordism Category associated with the Khovanov category $Cat(K).$ {\it We can investigate under what circumstances the Khovanov Cobordism Categores before and after the Reidemeister move are chain homotopy equivalent}. In Figure~\ref{diag} the maps labeled $H_{i}$ are homotopies. They are what you see. If a circle needs to be born or needs to be destroyed that is what the surface cobordism does.
In Figure~\ref{diag1} you see how the chain homotopy between $F_{1}G_{1}$ and the Identity is assembled from surface cobordisms. And you see that for the chain homotopy to satisfy
$\partial H + H \partial = F_{1}G_{1} + 1,$ the sum of mappings shown is needed. 

Now examine Figure~\ref{diag2} and you will see that the pattern of that requirement can be met if in the category of maps via 
surface cobordisms the so-called {\it 4-Tube Relation} is satisfied. In this figure, at the lower left, we see four local bits of surface labelled $1,2,3,4.$ This part of the figure is {\it not} part of the 4Tube Relation, but rather
an illustration of the four local surfaces that can be connected by tubes in the four ways.
This 4-Tube Relation says that if four bits of surface are nearby, call them $S_{1},S_{2},S_{3},S_{4}$ then you can form $S_{ij}$ by glueing a tube from 
$S_{i}$ to $S_{j}$ and then the relation is $$S_{12} + S_{34} = S_{14} + S_{23}.$$ In the figure we have shown tubes for $S_{ij}$ by constructing them between the corresponding surfaces.
When this tubing relation is satisfied then the chain homotopy class of the Bar-Natan Cobordism Category will be invariant under the 
Reidemeister Two move. One can prove that it will also be invariant under the Reidemeister Three move and  analyze the degree shift that results from the Reidemiester One move.
This then turns out to be enough to assemble the new invariant. 

Now an injunction for the reader: Examine Figure~\ref{diag1} and Figure~\ref{diag2}. Notice that Figure~\ref{diag2} has the {\it pattern} of the relation shown in Figure~\ref{diag1}. The pattern is not too complex.
Four bits of surface are near each other and tubes are constructed between pairs of the surfaces in four ways. One can do this construction with any collection of surfaces labeled $1,2,3,4.$ But in Figure~\ref{diag1}
we found this pattern by looking for a chain homotopy that would make our theory invariant under the second Reidemeister move. Our tubed surfaces came from assembling morphisms in the Khovanov Category of a 
knot or link. This part needs a lot of study and we have only introduced it with a sketch. We have done this to give you, the reader, a glimpse of how the general pattern of the 4Tube Relation is in back of the topological
invariance of Khovanov Homology. If you just see that this tubing relation is sufficient to give a formula at the Figure~\ref{diag1} level, that will be enough for a first pass. Then you can later read the references at the beginning of this section. There is a mystery here that the general pattern of the tubing relation is just what is needed to make the constructions work. 

It is worth remarking briefly how the invariance under the third Reidemeister moves comes about, as it is closely related to the way it happens for the original bracket polynomial. Let ${\cal C}(K)$ denote the Bar-Natan Cobordism Category of
a link diagram $K.$ Then we have a functor from the Bar-Natan Cobordism Category of an $A$-smoothing in $K$ to the Bar-Natan Cobordism Category of a $B$ smoothing in $K.$ 
$${\cal C}(\Across): {\cal C}(\Asmooth) \longrightarrow {\cal C}(\Bsmooth).$$
The functor takes arrows to arrows and objects to objects via the saddle point cobordism that re-smooths this one crossing.

\begin{figure}
     \begin{center}
     \begin{tabular}{c}
     \includegraphics[width=13cm]{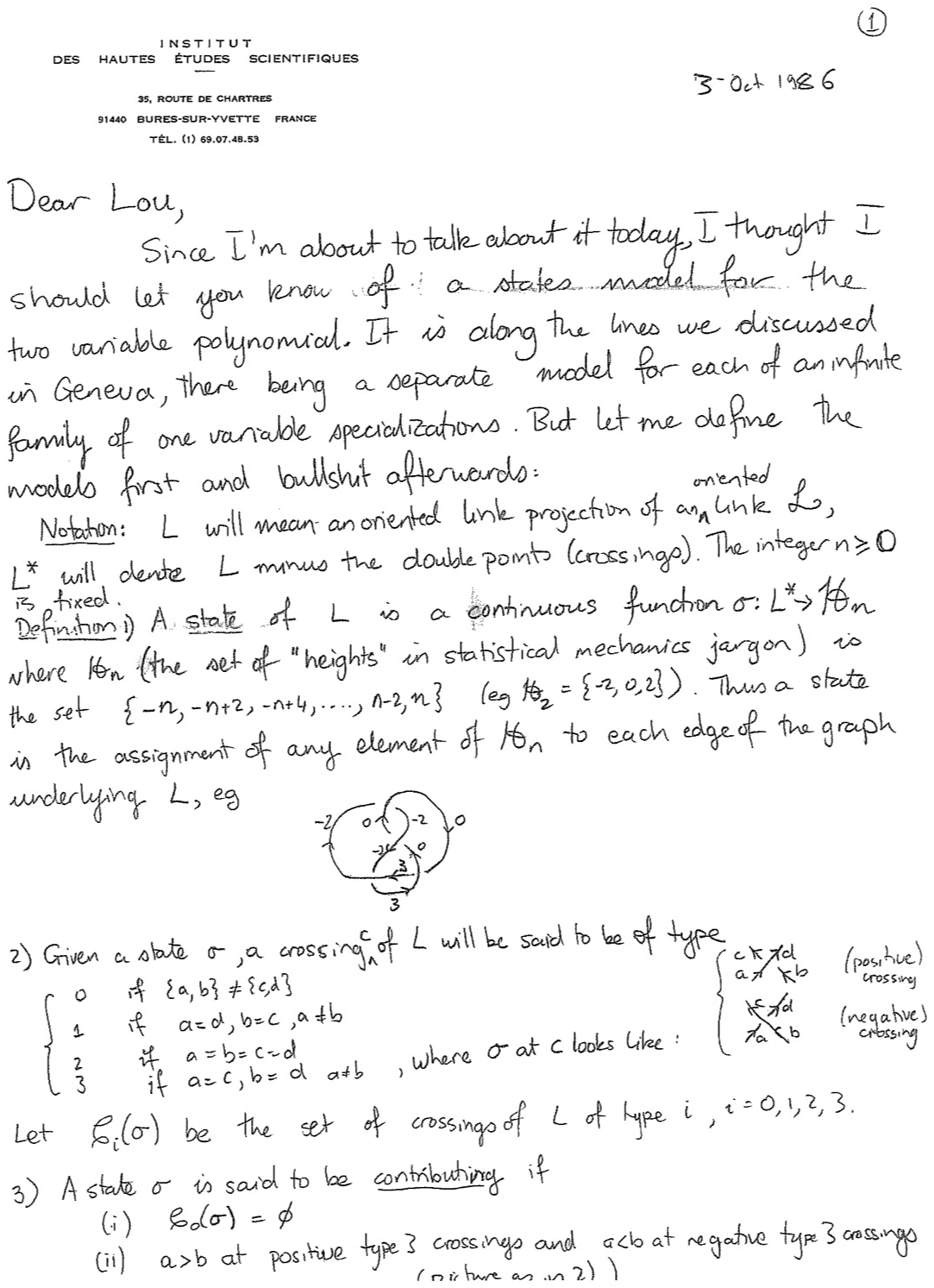}
     \end{tabular}
     \caption{\bf Vaughan Jones Letter}
    \label{VJ1}
\end{center}
\end{figure}

\begin{figure}
     \begin{center}
     \begin{tabular}{c}
     \includegraphics[width=13cm]{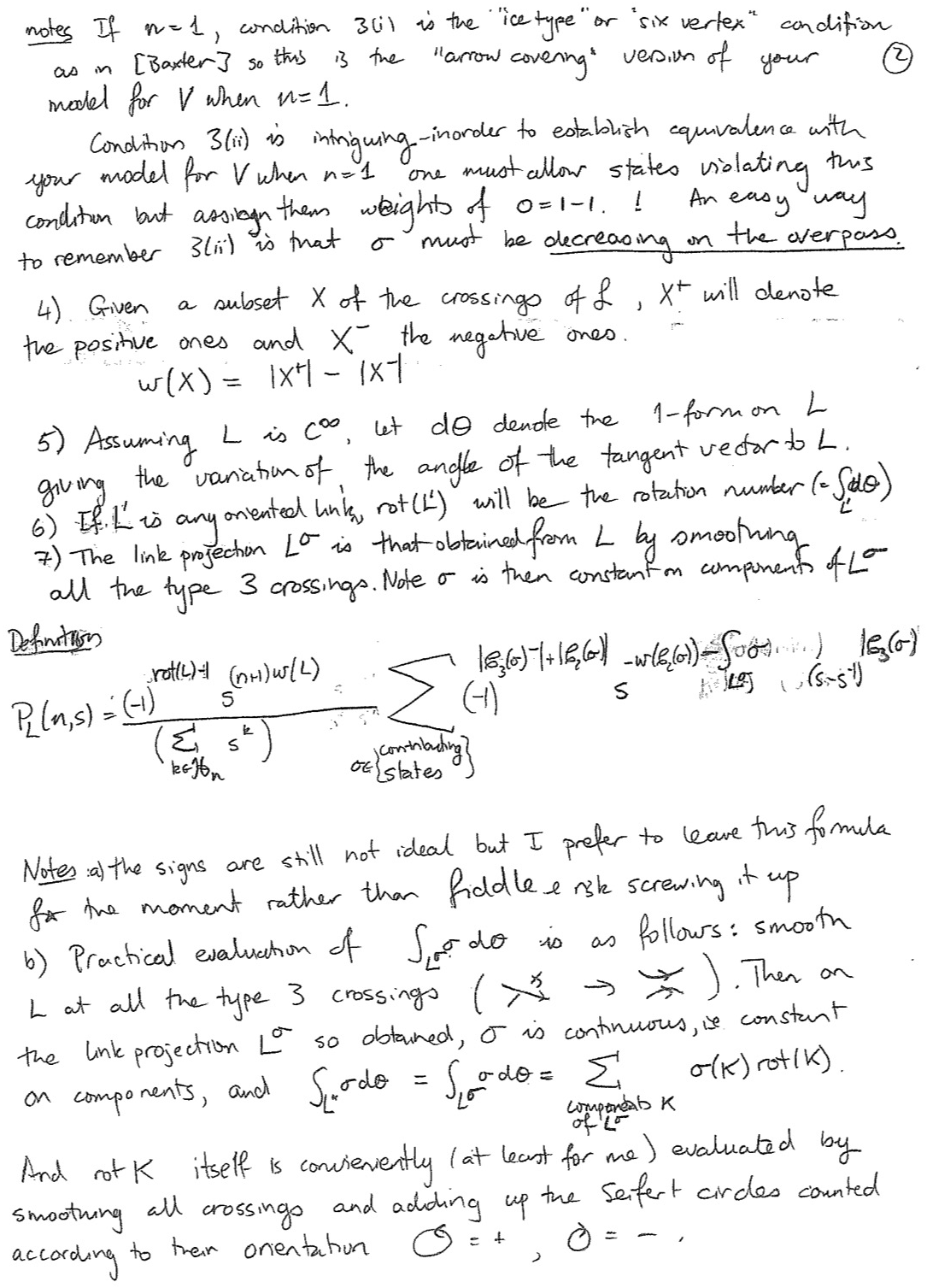}
     \end{tabular}
     \caption{\bf Vaughan Jones Letter}
    \label{VJ2}
\end{center}
\end{figure}

\begin{figure}
     \begin{center}
     \begin{tabular}{c}
     \includegraphics[width=13cm]{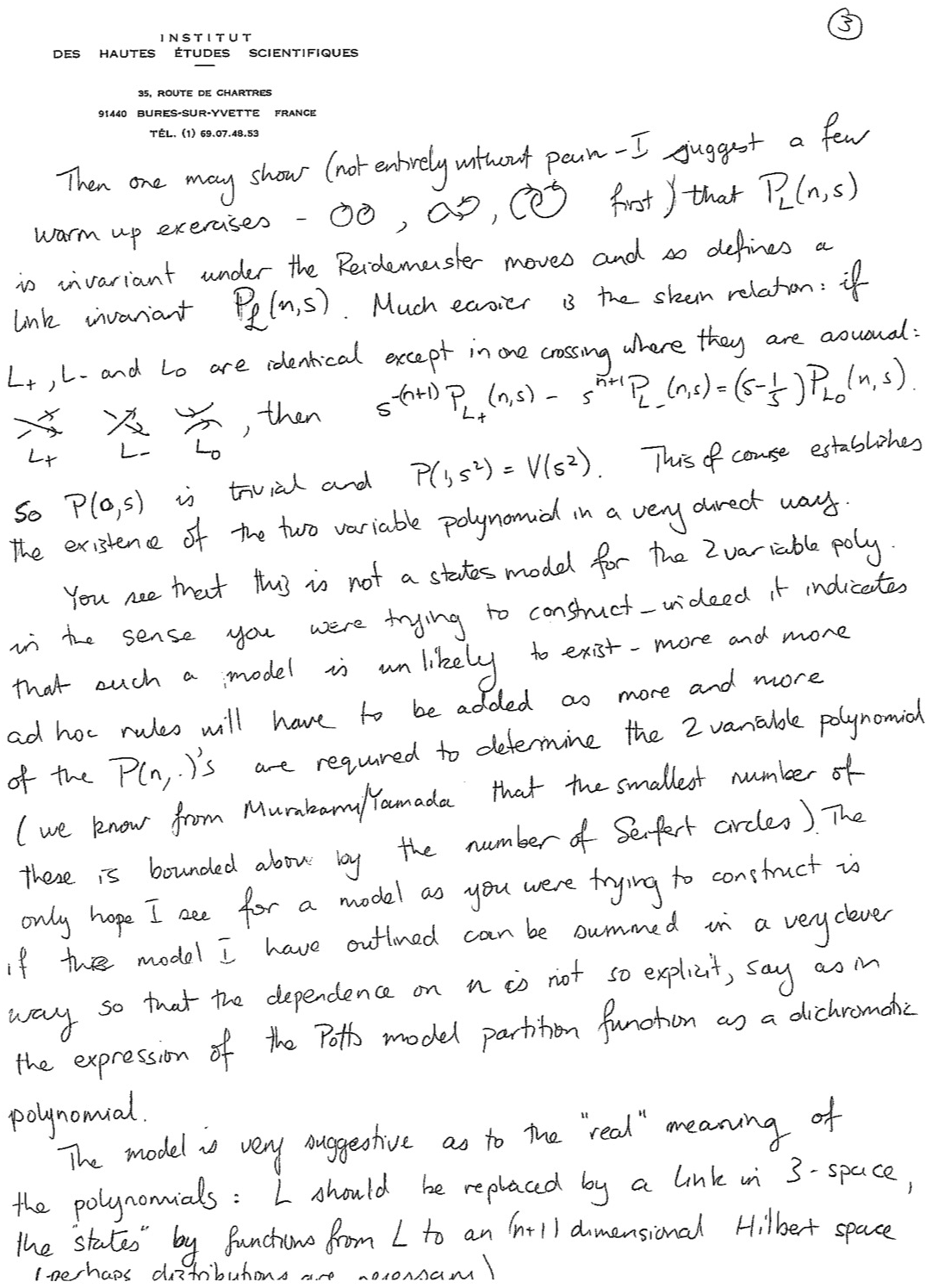}
     \end{tabular}
     \caption{\bf Vaughan Jones Letter}
    \label{VJ3}
\end{center}
\end{figure}

\begin{figure}
     \begin{center}
     \begin{tabular}{c}
     \includegraphics[width=13cm]{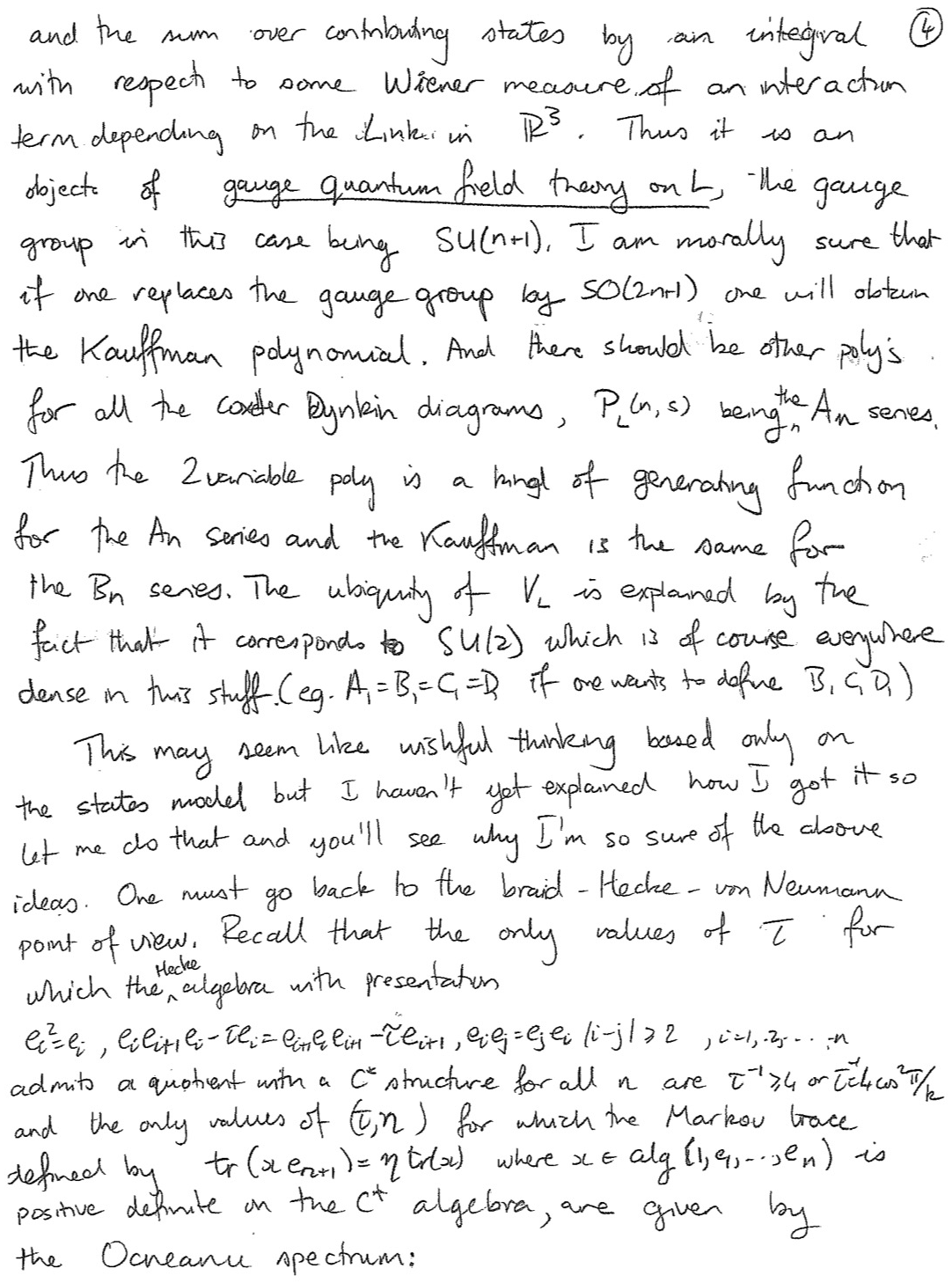}
     \end{tabular}
     \caption{\bf Vaughan Jones Letter}
    \label{VJ4}
\end{center}
\end{figure}

\begin{figure}
     \begin{center}
     \begin{tabular}{c}
     \includegraphics[width=13cm]{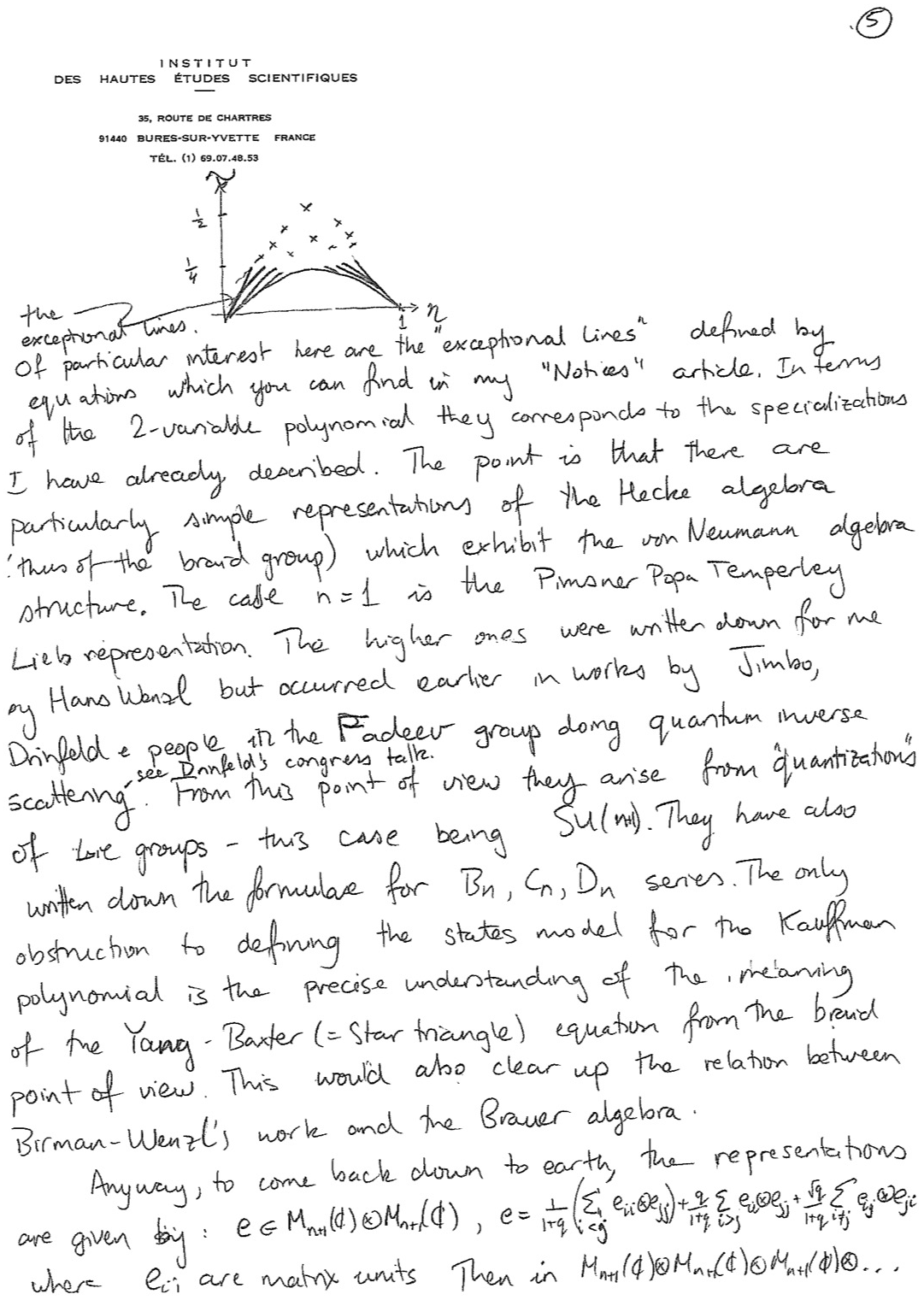}
     \end{tabular}
     \caption{\bf Vaughan Jones Letter}
    \label{VJ5}
\end{center}
\end{figure}

\begin{figure}
     \begin{center}
     \begin{tabular}{c}
     \includegraphics[width=13cm]{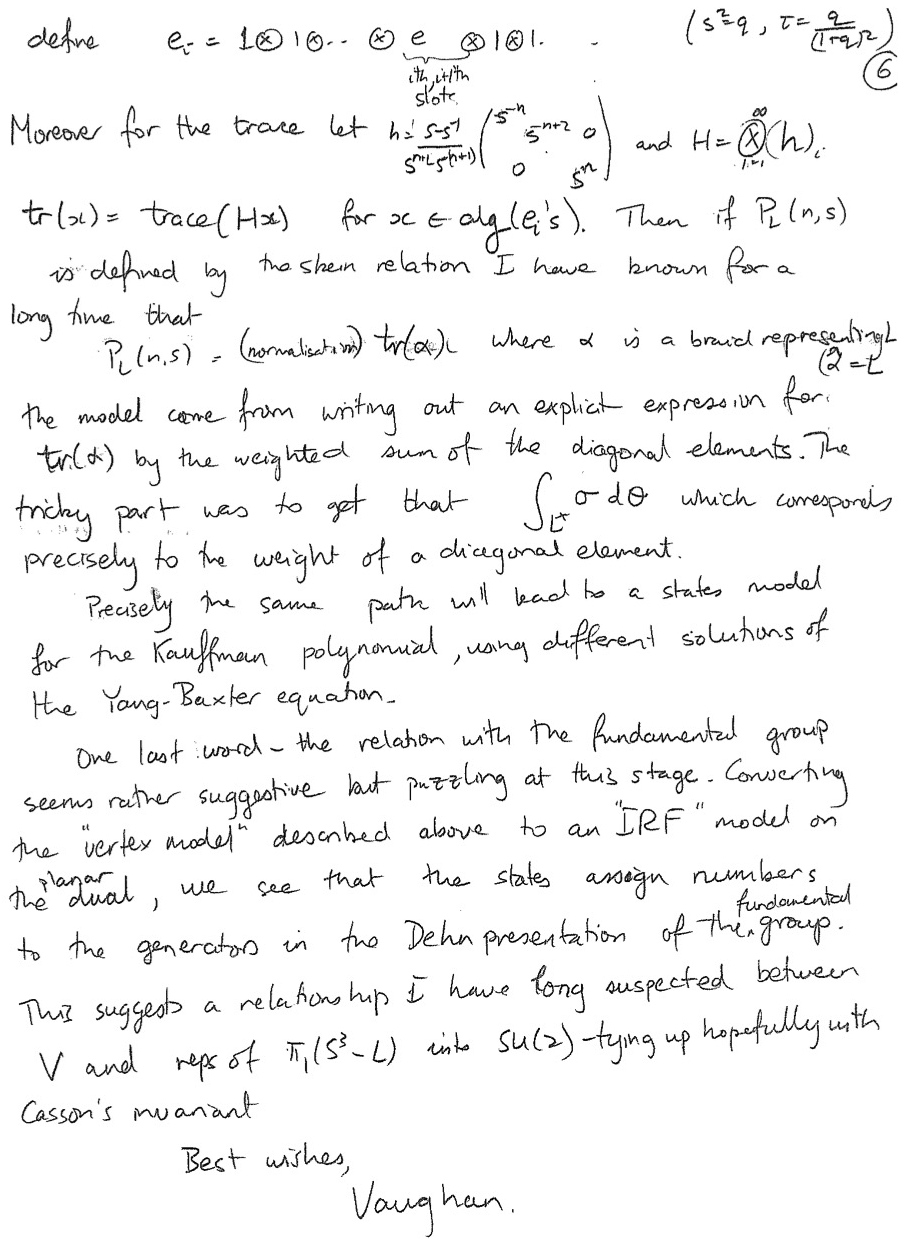}
     \end{tabular}
     \caption{\bf Vaughan Jones Letter}
    \label{VJ6}
\end{center}
\end{figure}

A functor from one category to another gives rise to a higher category since arrows in one category are taken to arrows in the other category by the ``arrow" of the functor. Since the functor takes objects to objects we
can flatten this higher category to a new category where an object in one category has an arrow from it to its functorial image in the other category. Thinking this way we see that 
the flattened higher category associated with the functor ${\cal C}(\Across): {\cal C}(\Asmooth) \longrightarrow {\cal C}(\Bsmooth)$ is exactly the Bar-Natan Cobordism Category ${\cal C}(K)$ where $K$ is the link whose crossing produced the functor. For example, consider the $4$-cube category in Figure~\ref{4cube}. In this figure we have drawn arrows on the edges from the outer $3$-cube to the inner $3$-cube, indicating a functor from outer $3$-cube to inner one. One can think of these arrows as the result of making a category from a functor as described above.\\

In fact, we can assemble the Bar-Natan Cobordism Category for ${\cal C}(K)$  from ${\cal C}(\Asmooth)$ and ${\cal C}(\Bsmooth)$ by forming the direct sum  ${\cal C}(\Asmooth) \bigoplus {\cal C}(\Bsmooth).$ More precisely, since the Bar-Natan Cobordism Category is graded by the number of $B$-smoothings in the states, we can write $${\cal C}^{n}(K) = {\cal C}^{n}(\Asmooth) \bigoplus {\cal C}^{n-1}(\Bsmooth).$$
If $\tau: {\cal C}^{n}(\Asmooth) \longrightarrow {\cal C}^{n-1}(\Bsmooth)$ denotes the resmoothing map, then we can write the boundary mapping on the direct sum by the formula
$$\partial (x,y) = (\partial x, \tau x + \partial y),$$ where for simplicity we have written this formula as though there were elements and we have written it modulo two. This formula expresses the fact that the boundary for 
a given state is obtained either by locally smoothing the state at the indicated site (this is $\tau$) or by using smoothings elsewhere in the state that are not directly indicated (this is the $\partial$ part of the formula).

One recognizes that this is the 
familiar mapping cone construction from homological algebra. Now view Figure~\ref{mcone}. Here we indicate that there is a functor from the Bar-Natan Cobordism Category of $K$, $A$-smoothed at a Reidemeister three move configuration, to the corresponding resmoothing. Reidemeister $2$-moves map the source Bar-Natan Cobordism Category to a chain-homotopy equivalent one. One can verify that the corresponding direct sum Khovanov Cobordism Categories are chain homotopy equivalent. From this it follows that the chain homotopy class of the Bar-Natan Cobordism Category of a link diagram $K$ is not changed by the third Reidemeister move. This proof is a categorified version of the original proof that the 
bracket polynomial is invariant under regular isotopy.

Bar-Natan Cobordism Category is the categorical structure behind Khovanov Homology. Functors from the Bar-Natan Cobordism Category to categories of modules can be constructed and actual homology calculated. 
The Bar-Natan Cobordism Category gives a diagrammatic/categorical understanding of how this homology theory gives topological information about knots and links.

It is not lost on us that the 4-Tube Relation has an analogy with the the 4-Term Relation in the theory of Vassiliev invariants (as we have described it in Section 3). The 4-Term Relation is closely tied with the Jacobi identity in Lie algebras, and specific Lie algebras can be used to construct invariants of knots and links. If one follows this analogy with the 4-Tube Relation, as Dror Bar-Natan did, one finds that there are certain
Frobenius algebras (see \cite{LKKho,BN2,DKK}) that are instrumental for constructing link homology theories. The original functor devised by Khovanov can be described very simply in terms of such a Frobenius algebra.
Let $V = Z[x]/(x^2)$ be the polynomial ring over the integers with transcendental variable $x$ modulo the ideal generated by $x^2.$
If a state of the bracket polynomial has $k$ loops send it to the $k$-fold tensor product of $V$ with itself. The morphisms of the Khovanov category of a knot K are, as we know, described by resmoothings at crossings.
Such a resmoothing corresponds to a surface cobordism taking two loops to one loop or to a surface cobordism taking one loop to two loops. These morphisms are illustrated in Figure~\ref{saddle1}. Let 
$m$ denote the morphism from two loops to one loop and $\Delta$ denote the morphism from one loop to two loops.  Then the functor that takes loops to tensor products of the algebra $V$ will take $m$ to 
the multiplication in the algebra and $\Delta$ to the operation defined below that is a co-multiplication on the algebra. We use the symbol $\Delta$ again for its image under the functor.
\begin{enumerate}
	\item $\Delta(1) = 1 \otimes x + x \otimes 1.$
	\item $\Delta(x) = x \otimes x.$
\end{enumerate}
It is not hard to verify that this indeed defines a functor from the Khovanov Category of a knot $K$ to a category of modules over the integers and that the functor is compatible with theBar-Natan Cobordism Category construction,
so that the image of the Bar-Natan Cobordism Category of $K$ is a chain complex. The homology of this chain complex is Khovanov homology.
We will not here go into further algebraic details. The Jones polynomial itself is a graded Euler characteristic of the 
Khovanov homology. This part of the development is quite analogous to the original development of homology theory where the homology groups replaced the Betti numbers. The Khovanov homology is, after all is said and done, a natural categorification of the Jones polynomial.

For the cobordism point of view that we have discussed here, we recommend the paper by Bar-Natan \cite{BN2}.

Through its categorification, the Jones polynomial reaches a great result. Kronheimer and Mrowka \cite{KM} proved that the Khovanov Homology detects knottedness of classical knots.
It remains an open problem at this writing whether the orginal Jones polynomial detects knottedness. Is there a non-trivial classical knot $K$ with unit Jones polynomial?

The Khovanov homology is one example of a number of {\it link homology} theories that have been discovered. At this point we can  mention the Heegaard-Floer Link Homology \cite{Ozsvath1} whose chain complex can be formulated in terms of the spanning tree states with which we began this essay. The Heegaard-Floer Link Homology is originally defined by Ozsvath and Szabo in terms of high dimensional symplectic
geometry and a chain complex whose generators are in 1-1 correspondence the the Formal Knot Theory \cite{FKT} states for the Alexander Conway polynomial. The differentials in the complex originally did not have 
a definition in terms of the FKT states. This has been rectified in recent work such as \cite{Ozsvath2,Baldwin} .

In all of this development the Jones polynomial has been the keystone in guiding researchers to the right 
constructions and the most creative questions.

\section{The Vaughan Jones Letter of 1986}
This letter of Vaughan Jones is significant in that it encapsulates Vaughn Jones' extraordinary vision of the possiblities for these invariants.
Papers by Jones and by Turaev and Wadati soon appeared, showing how to generate many knot invariants, including the original Jones polynomial by using solutions of the Yang-Baxter Equation.
This cemented the direct relationship with statistical mechanics models that had begun with the state summation model and the Kauffman bracket. The letter indicates in detail the persistence with which
Jones followed this subject of his own making. The complete letter is given in Figures~\ref{VJ1},\ref{VJ2},\ref{VJ3},\ref{VJ4},\ref{VJ5},\ref{VJ6}.\\

\end{document}